\documentclass[review]{elsarticle}
\usepackage[utf8]{inputenc}

\usepackage{orcidlink}
\usepackage{amsmath,amssymb}
\usepackage{bm}
\usepackage{multirow}
\usepackage{algpseudocode}

\usepackage{hyperref}
\usepackage[misc]{ifsym} 

\usepackage{color}
\usepackage{graphicx}
\graphicspath{{./figures/}}
\usepackage{subcaption}

\usepackage{amsfonts}
\usepackage{mathtools}

\usepackage[mathscr]{eucal}

\usepackage{pgfplots}

\usepackage{bbm}
\ifpdf
\DeclareGraphicsExtensions{.eps,.pdf,.png,.jpg}
\else
\DeclareGraphicsExtensions{.eps}
\fi

\newtheorem{remark}{Remark}
\newtheorem{assumption}{Assumption}
\newcommand{\mat}[1]{\mathbf{#1}}

\providecommand{\rset}{\mathbb{R}} 

\providecommand{\E}[1]{{\ensuremath{\mathbb{E}}\mspace{-2mu}\left[#1\right]}}    
\providecommand{\var}[1]{{\ensuremath{\mathrm{Var}}\mspace{-2mu}\left[#1\right]}}

\providecommand{\QoI}{\ensuremath{Q}} 
\providecommand{\tol}{\ensuremath{\mathrm{TOL}}} 
\providecommand{\work}{\ensuremath{W}}           

\DeclarePairedDelimiter\abs{\lvert}{\rvert}%
\DeclarePairedDelimiter\norm{\lVert}{\rVert}%

\newcommand{\RN}[1]{%
	\textup{\uppercase\expandafter{\romannumeral#1}}%
}

\definecolor{revR1color}{rgb}{0.00,0.30,0.56}   
\definecolor{revAEcolor}{rgb}{0.09,0.48,0.25}   
\definecolor{revR2color}{rgb}{0.66,0.14,0.14}   
\definecolor{revAllcolor}{rgb}{0.50,0.00,0.50}  
\newcommand{\revRone}[1]{{#1}}
\newcommand{\revRtwo}[1]{{#1}}
\newcommand{\revAE}[1]{{#1}}

\providecommand{\workExp}{\ensuremath{\chi}}

\begin{document}

\begin{abstract}
  The efficient approximation of quantities of interest derived from PDEs with lognormal diffusivity is a central challenge in uncertainty quantification. \revRtwo{This paper targets a problem class that combines four analytical difficulties: a geometric boundary singularity, a lognormal coefficient field without a deterministic positive lower bound, sample-dependent mesh selection that introduces parameter-space discontinuities, and infinitely many discontinuity locations that preclude classical pre-integration smoothing.} In this study, we propose a multilevel quasi-Monte Carlo framework to approximate deterministic, real-valued, bounded linear functionals that depend on the solution of a linear elliptic PDE with a lognormal diffusivity coefficient parameterized by a \revAE{multi-dimensional} Gaussian random vector and deterministic geometric singularities in bounded domains of $\mathbb{R}^d$. We analyze the parametric regularity and develop the multilevel implementation based on a sequence of adaptive meshes, developed in \revAE{our earlier work} ``Goal-oriented adaptive finite element multilevel Monte Carlo with convergence rates'', CMAME, 402 (2022), p. 115582. For further variance reduction, we incorporate importance sampling and introduce a level-0 control variate within the multilevel hierarchy. \revAE{Introducing such a control variate can alter the optimal choice for the initial mesh, further highlighting the advantages of adaptive meshes.} \revAE{On a 2-D slit benchmark discretized with bilinear, quadrilateral Q1-FEM, numerical experiments show that, in the parameter range explored, the proposed adaptive MLQMC algorithm achieves a prescribed accuracy at markedly lower computational cost than a standard multilevel Monte Carlo estimator on the same mesh hierarchy.}
\end{abstract}

\title{Goal-Oriented Adaptive Finite Element\\ Multilevel Quasi-{M}onte {C}arlo}

\author[KFUPM]{Joakim~Beck\,\orcidlink{0000-0001-8042-998X}}
\author[CEMSE]{Yang~Liu\corref{cor1}\fnref{fn2}\,\orcidlink{0000-0003-0778-3872}}
\ead{yang.liu.3@kaust.edu.sa}
\author[CEMSE]{Erik~von~Schwerin\fnref{fn2}\,\orcidlink{0000-0002-2964-7225}}
\author[CEMSE]{Ra\'{u}l~Tempone\fnref{fn2,fn3}\,\orcidlink{0000-0003-1967-4446}}
\address[KFUPM]{College of Petroleum Engineering \& Geosciences, Center for Integrative Petroleum Research, King Fahd University of Petroleum and Minerals, Dhahran 31261, Kingdom of Saudi Arabia}
\address[CEMSE]{Computer, Electrical and Mathematical Sciences and Engineering,\\
  4700 King Abdullah University of Science and Technology (KAUST),\\
  Thuwal 23955-6900, Kingdom of Saudi Arabia}

\cortext[cor1]{Corresponding author}

\fntext[fn2]{KAUST SRI Center for Uncertainty Quantification in Computational Science and Engineering}
\fntext[fn3]{Alexander von Humboldt Professor in Mathematics for Uncertainty
  Quantification, RWTH Aachen University, 52062 Aachen, Germany.}

\markboth{Beck et al.}{Goal-Oriented Adaptive Multilevel Quasi {M}onte {C}arlo}

\maketitle

\section{Introduction}
\label{sec:intro}

We consider a physical system subject to uncertainty and modeled by a random partial differential equation (RPDE) 
\cite{Cetal2017,OBF2017,lord_powell_shardlow_2014}, and a given scalar quantity of interest (QoI), depending 
on the solution to a boundary value problem for the RPDE.
The topic of this work is adaptive computation and error control for QoI expectations of the form $\E{\QoI(u)}$, where 
$\QoI$ is a deterministic, real-valued, bounded linear functional of $u$ which almost surely solves the boundary value problem 
of a linear elliptic partial differential equation (PDE) with random coefficients:
\begin{subequations}
  \label{eq:bvp_general}
  \begin{align}
    -\nabla \cdot \left( a(\mathbf x; \omega) \nabla u(\mathbf x; \omega) \right) 
    &= f(\mathbf x) &&\text{for $\mathbf x  \in \mathcal D$,} 
    \\
    u(\mathbf x; \omega) &=  0 &&\text{for $\mathbf x \in \partial\mathcal{D}_1$},
    \\
    \partial_n u(\mathbf x; \omega) &=  0 &&\text{for $\mathbf x \in \partial\mathcal{D} -\partial\mathcal{D}_1$}.
  \end{align}
\end{subequations}
The variable $\omega$ corresponds to an outcome associated with a complete probability space $(\Omega,\mathcal{F},\mathbb P)$
and the variable $\mathbf x$ belongs to an open and bounded polygonal/polyhedral domain $\mathcal{D}$ in 
$\rset^d$, where $d\ge2$. The divergence and gradient operators $\nabla\cdot$ and $\nabla$ are applied with respect to the spatial variable $\mathbf x$.
The randomness in the stochastic diffusivity coefficient field $a(\mathbf x;\omega)$ in general causes 
the solution $u$ to be stochastic.

The boundary $\partial\mathcal{D}$ is divided into two disjoint parts with homogeneous Dirichlet and Neumann boundary conditions, respectively. Both parts of the boundary are unions of a finite number of intervals or polygons. \revRone{This choice reflects the physical setup of the needle problem as in~\cite{adMLMC_our}. A geometric singularity arises at the fixed, a priori known interface between the Dirichlet and Neumann regions; this singularity drives our adaptive mesh refinement to concentrate elements near that interface.}


An important class of RPDEs in applications are those with lognormal coefficient fields, representing
physical properties such as conductivity, diffusivity, or elasticity. 
This class is a typical choice for these 
fields because the fields are strictly positive and the distribution has a heavy tail, often observed in practice; see~\cite{LLSY2019}. 
This work focuses on lognormal coefficient fields parametrized by a finite number of random variables, 
of the form
\begin{align}
  \label{eq:a_def}
  a(\mathbf{x};\mathbf{y}) = \exp{\left({y}_1 + \sum_{j = 2}^{s} {y}_j \psi_j(\mathbf{x})\right)},
\end{align}
where $s \in \mathbb{N}$ and $\mathbf{y}_j$, for $j = 1,\dotsc,s$, are independent 
standard normal random variables. 
The functions $\psi_n$, $n=1,\dots,s$, correspond to a series representation of the random field, which we 
assume to be sufficiently differentiable for the adaptive finite element method (FEM) error estimator 
introduced in previous work~\cite{adFEM_our}, with a correlation length comparable to the domain size. 
\revAE{A lognormal coefficient field is strictly positive almost surely, but it admits no deterministic positive lower or upper bound. As a consequence, the bilinear form associated with~\eqref{eq:bvp_general} is coercive and continuous pathwise, but the constants of coercivity and continuity are not uniform in $\omega$; the well-posedness analysis must therefore proceed through $\omega$-dependent constants with finite moments, as carried out in detail in~\cite{BNT_SIAM_review,Charrier_2012,CD2013,Cliffe_2011,Teckentrup_2013}. This weakened coercivity propagates into the QoI in two ways that drive the present work: (i) the mapping $\mathbf y \mapsto Q$ can become unbounded along low-probability directions, leading to heavy-tailed sample variance for Monte Carlo; and (ii) standard QMC theory based on the unweighted Koksma--Hlawka inequality no longer applies, because the integrand is no longer in the unit-cube Sobolev class assumed by the classical bounds.}

We could also consider a stochastic forcing, $f(\mathbf x; \omega)$, as long as we assume that any 
need for highly localized adaptive mesh refinement is driven by the deterministic geometry of $\mathcal{D}$.
In this setting it is favorable to generate deterministic $h$-adaptively refined meshes, 
adapted to the geometrically induced singularity and the QoI, and use sample-adaptive 
selection of such meshes as described in~\cite{adMLMC_our}.

The multilevel Monte Carlo (MLMC) method employs a hierarchy of discretizations to reduce the complexity of traditional Monte Carlo (MC) simulations. This approach was initially introduced independently by Heinrich~\cite{heinrich2001multilevel} and Giles~\cite{Giles_OpRes} in different contexts. Giles' approach extended prior work of Kebaier~\cite{kebaier05}, who applied two-level discretizations of stochastic differential equations (SDEs) as control variates to reduce computational complexity. A comprehensive overview of MLMC methodologies is available in~\cite{AcNum_MLMC}. Some adaptations of the multilevel hierarchy include the optimization of MLMC strategies~\cite{optimal_hierarchies_our} and the development of the Continuation MLMC approach~\cite{continuation_MLMC_our}. 

MLMC has been successfully applied to numerical approximations of RPDEs, especially elliptic PDEs with random coefficients, as shown in~\cite{adMLMC_our, MLMC_FEM_Schwab,Charrier_2012,CD2013,cst13,Cliffe_2011,Teckentrup_2013}. \revAE{The work~\cite{Hall_Hoel} developed computable error estimates for elliptic PDEs with rough stochastic data in the single level setting.}

Multilevel Quasi-Monte Carlo (MLQMC) methods replace MC sampling with deterministic quadratures whose low discrepancy yields faster convergence for each level estimator. For integrands with sufficient regularity, randomized QMC points reduce the variance more efficiently than MC at each level, improving overall computational complexity. MLQMC is first proposed for stochastic differential equations simulations~\cite{giles2009multilevel} and MLQMC for RPDEs with uniform mesh refinements have been proposed in~\cite{HS2019,KSSSU2017}, while the analysis of QMC in this setting is treated in~\cite{Graham_et_al_Num_Math_2015}. 

\revRtwo{A complementary line of work approaches the same setting of elliptic RPDEs through stochastic collocation on structured grids, where the integrand is tabulated on a tensorized or sparse set of collocation nodes and either integrated or approximated~\cite{BNT_SIAM_review,NobileTemponeWebster2008b,NobileTemponeWebster2008a}. Anisotropic and multilevel variants extend this idea to problems whose mixed regularity is uneven across parameter directions. A multilevel stochastic collocation analysis with convergence and cost estimates for elliptic PDEs, where the coefficients field is uniformly bounded away from 0~\cite{TeckentrupJantschWebsterGunzburger2015}. These methods are most competitive when the parametric integrand has sufficient mixed regularity and its effective dimensionality is moderate; in the present paper the combination of a geometric boundary singularity, a lognormal (unbounded) diffusion, and sample-dependent $h$-adaptive meshes produces low parametric regularity and discontinuities in $\mathbf y \mapsto Q$, regimes in which sparse-grid collocation loses its asymptotic advantage. A quantitative comparison on the present slit benchmark is outside the scope of this paper.}

In RPDEs with lognormal coefficients, the mapping $\mathbf{y} \mapsto Q$ from random parameters to the QoI exhibits boundary singularities, and the classical Koksma--Hlawka inequality cannot be applied for the error estimation. In this situation, weighted Sobolev spaces are used in~\cite{Graham_et_al_Num_Math_2015} to address the singularity, with corresponding lattice rules designed accordingly. Additionally, Owen~\cite{Owen06} introduced the boundary growth condition to characterize boundary singularities. Recent work on RPDEs with lognormal coefficients in~\cite{liu2023nonasymptotic} has identified a QMC convergence rate of $\mathcal{O}(N^{-1+\epsilon})$. Further investigations into RQMC convergence rates for integrands with boundary unboundedness and interior discontinuities can be found in~\cite{liu2025weak,liu2024randomized}. 

As discussed in~\cite{adMLMC_our}, stochastic mesh selection is highly beneficial for MLMC when applied to~\eqref{eq:bvp_general} with the lognormal coefficient field~\eqref{eq:a_def}. However, sample-dependent mesh selection introduces discontinuities in the mapping $\mathbf{y} \mapsto \QoI$, as shown in~\cite{adMLMC_our}. While these discontinuities do not affect the complexity rate for MC-based methods, they reduce the benefits of QMC-based methods. Some previous work, e.g.~\cite{griewank2018high}, address discontinuities by employing a pre-integration smoothing method for integrands of the form $f \mathbbm{1}_g$, where $\mathbbm{1}_g$ is an indicator function that equals 1 when $g > 0$, and $0$ otherwise. A key assumption in this approach is the strict monotonicity of $g$ with respect to certain variables. The works~\cite{Beyer_Ben_Hammouda_Tempone_2020,Bayer_Siebenmorgen_Tempone_2018} are established on this setting. However, as analyzed in~\cite{gilbert2022preintegration}, the pre-integration smoothing fails to yield a function within the desired Sobolev space when the monotonicity condition is not satisfied. 

In our study, the discontinuous integrand $Q$ does not satisfy the monotonicity assumption. Moreover, the number of discontinuities is infinite \revAE{due to the lognormal coefficient field} and their locations are unknown a priori. If discontinuity locations were known, one could integrate over each continuous region and aggregate the results. We examine the feasibility of this approach in~\ref{sec:adaptivity_qmc_piecewise_integration}.

\revRtwo{The problem class targeted in this paper simultaneously exhibits four analytical difficulties that, when combined, preclude the direct application of the standard MLMC or MLQMC complexity theory cited above:
\begin{itemize}
\item[(D1)] \emph{Geometric singularity at a fixed interface.} The slit/needle geometry of $\partial\mathcal{D}_1$ produces a deterministic corner singularity in $u$ whose strength is independent of $\omega$; this demands strongly graded or $h$-adaptive meshes near the interface and makes convergence using uniform meshes suboptimal.
\item[(D2)] \emph{Lognormal (unbounded) coefficients.} As discussed above, the bilinear form is not uniformly coercive and the QoI can be heavy-tailed; all error and complexity bounds must tolerate coefficients without a deterministic lower bound.
\item[(D3)] \emph{Sample-dependent stopping in the adaptive FEM.} The pathwise mesh $\mathcal{T}_\ell(\mathbf y,\tol_\ell)$ depends on the realization through the a posteriori error estimator, so the map $\mathbf y \mapsto Q$, from parameters to QoI, is piecewise smooth with discontinuities across the level surfaces of the stopping criterion; see Figure~\ref{fig:mesh_selection_discontinuity}.
\item[(D4)] \emph{Infinitely many unknown discontinuity locations.} The union of those level surfaces is dense enough that classical pre-integration smoothing~\cite{griewank2018high} does not apply, because no single coordinate direction is known a priori to be a monotonicity direction for the indicator.
\end{itemize}
Difficulty (D1) is addressed by goal-oriented $h$-adaptivity~\cite{adFEM_our}; difficulty (D2) by working with moment based well-posedness estimates of~\cite{Charrier_2012,CD2013,Cliffe_2011,Teckentrup_2013}; difficulties (D3) and (D4) are the focus of Sections~\ref{sec:MC_QMC_estimators}--\ref{sec:variance_reduction} of the present paper.}

\revRtwo{\paragraph{Contributions} With difficulties (D1)--(D4) in mind, the main contributions of this paper are the following:
\begin{itemize}
\item[(C1)] \emph{An MLQMC method for elliptic PDEs with lognormal coefficients and a geometric boundary singularity, built on goal-oriented $h$-adaptivity.} We combine the sample-adaptive mesh-selection strategy of~\cite{adMLMC_our} with randomized-shift QMC quadrature at each level and analyse the resulting computational complexity in the three regimes dictated by the signs of $\beta-\workExp$ and $\beta-\lambda\,\workExp$ (see Table~\ref{tab:complexity_summary} and Assumption~\ref{assumption:qmc_rates}), where $\workExp \coloneqq \gamma d$ is the effective work-growth exponent introduced in Assumption~\ref{assumption:mc_complexity}. In the 2-D slit problem with bilinear ($Q_1$) Lagrange FEM on quadrilateral meshes, adaptive meshes reduce $\workExp$ to 1 from 2 with uniform meshes, consistent with the work per level scalings derived in Section~\ref{sec:MC_QMC_estimators}. The observed MLQMC cost lies between $\mathcal{O}(\tol^{-1})$ and $\mathcal{O}(\tol^{-2})$, quantified in~\eqref{eq:mlqmc_adaptive_2d_complexity}.
\item[(C2)] \emph{Variance reduction tailored to the discontinuous integrand.} We design importance sampling proposals~$\varphi_{\bm{\alpha}}$ that push probability mass away from the boundary singularities of the integrand, and we combine them with two families of level $0$ control variates (truncation-based with importance-ranked or natural-order selection, and SVD-based). 
\item[(C3)] \emph{Quantitative numerical evidence on a 2-D slit benchmark.} Section~\ref{sec:numex} compares adaptive MLQMC against adaptive MLMC, non-adaptive MLQMC, and the variance-reduction variants on the 2-D slit benchmark, and reports the observed cost-versus-tolerance curves, the variance decay across levels, and the effect of importance sampling and control variates. The results confirm the theoretical regimes of Table~\ref{tab:complexity_summary} and expose the parameter ranges in which each variance-reduction building block is most effective.
\end{itemize}}

\revAE{The numerical study is centered on a representative 2-D slit benchmark, where the mixed Dirichlet--Neumann boundary conditions give rise to a fixed geometric singularity. The lognormal diffusivity is represented by a truncated expansion based on a Matérn covariance model and driven by a 49-dimensional standard Gaussian vector. This benchmark provides a concrete test case in which the analytical difficulties described above occur together in a controlled numerical setting.}

The remainder of this work is organized as follows: Section~\ref{sec:prob} formulates the problem and introduces the adaptive MLQMC method. Section~\ref{sec:MC_QMC_estimators} presents the MLMC and MLQMC estimators and analyzes their computational complexity. Section~\ref{sec:variance_reduction} describes variance reduction techniques, including importance sampling and control variates. Section~\ref{sec:numex} presents the numerical results and Section~\ref{sec:concl} concludes this study.


\section{Problem setting}
\label{sec:prob}
We recall the random elliptic PDE model~\eqref{eq:bvp_general},
\begin{subequations}
	\label{eq:bvp_general_1}
	\begin{align}
	\label{eq:pde_general_1}
	-\nabla \cdot \left( a(\mathbf x; \mathbf y) \nabla u(\mathbf x; \mathbf y) \right) 
	&= f(\mathbf x) &&\text{for $\mathbf x  \in \mathcal D$,} 
	\\
	u(\mathbf x; \mathbf y) &=  0 &&\text{for $\mathbf x \in \partial\mathcal{D}_1$},
	\\
	\partial_n u(\mathbf x; \mathbf y) &=  0 &&\text{for $\mathbf x \in \partial\mathcal{D} -\partial\mathcal{D}_1$},
	\end{align}
\end{subequations}
where the randomness is \revAE{represented by the random vector} $\mathbf{y}$. The mixed boundary setting is particularly relevant in applications such as the slit problem, where sharp geometric features lead to singularities in the solution and motivate the need for adaptive finite element methods (FEM) to resolve these features. 
In this work, $a(\mathbf{x};\mathbf y)$ takes the form
\begin{align*}
	a(\mathbf{x};\mathbf y) = \exp\left({y}_1  + \sum_{j = 2}^{s} {y}_j \psi_j (\mathbf{x})\right),
\end{align*}
where $s \in \mathbb{N} $, $\mathbf y = (y_1, \dots, y_s)$, and $y_j$ are independent random variables with the distribution 
$\mathcal{N}(0,1)$, for $j = 1,\dots, s$. 
We assume the basis functions $\psi_j$ to be sufficiently differentiable for the error estimator in adaptive FEM introduced in previous work~\cite{adFEM_our}. 
The random field also has a correlation length comparable to the domain size. Given $a$, as defined in~\eqref{eq:a_def}, 
the linearity of $Q$ in combination with~\eqref{eq:bvp_general_1} implies the following:
\begin{equation}
	Q(u(\mathbf{y})) = \exp(-y_1) Q(u(0;\mathbf{y}_{-1})),
\end{equation}
where $\mathbf{y}_{-1} \coloneqq (y_2,\dots,y_s)$. 
For simplicity, the notation $\tilde{Q}(\mathbf{y}_{-1}) \coloneqq Q(u(0; \mathbf{y}_{-1}))$ is used in the rest of this work. 
The aim is to compute the expectation of the QoI:
\begin{align}
	\mathbb{E}[Q(u(\mathbf{y}))] = \int_{\mathbb{R}^s} Q(u(\mathbf{y})) \varphi (\mathbf {y}) d\mathbf{y},
	\end{align}
where $\varphi (\mathbf y)$ denotes the $s$-dimensional standard normal distribution density.
Since $\mathbb{E}[\exp(-y_1)]$ is known, computing the expectation $\E{Q}$ reduces to computing $\E{\tilde{Q}(\mathbf{y}_{-1})}$. \revAE{From this point onward we adopt the standing convention that every multilevel estimator targets $\E{\tilde{Q}(\mathbf{y}_{-1})}$, and the expectation of the original QoI is recovered through $\mathbb{E}[Q(u(\mathbf y))] = \mathbb{E}[\exp(-y_1)]\,\mathbb{E}[\tilde{Q}(\mathbf{y}_{-1})] = e^{1/2}\,\mathbb{E}[\tilde{Q}(\mathbf{y}_{-1})]$. Throughout the rest of the paper, we therefore work exclusively with $\tilde{Q}$ and its level counterparts $\tilde{Q}_k$ and $\Delta\tilde{Q}_k$; the untilded symbol $Q$ is retained only in~\ref{sec:adaptivity_qmc_piecewise_integration}, where the discontinuity argument is made on the original QoI.}

\begin{figure}
\centering
\begin{subfigure}{0.47\linewidth}
	\includegraphics[width=\textwidth]{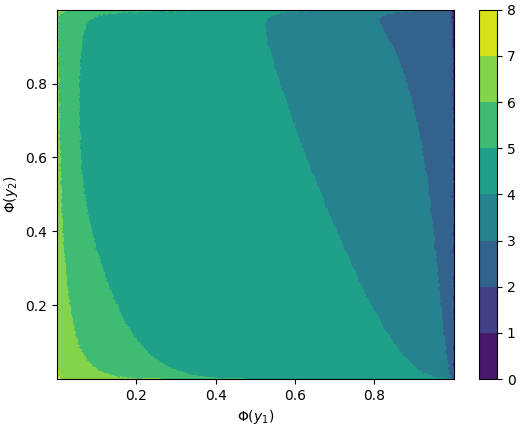}
	\caption{\revRtwo{Mesh index $\mathcal{K}(\mathbf{y}, \tol_0)$ on the two-parameter domain with the coefficient}}
\end{subfigure}
\begin{subfigure}{0.47\linewidth}
	\includegraphics[width=\textwidth]{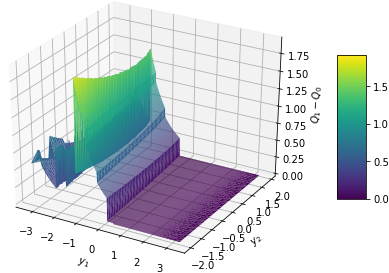}
	\caption{\revRtwo{Hierarchical difference $\abs*{{Q}(\mathbf{y}; \tol_1) - {Q}(\mathbf{y}; \tol_0)}$ in the parametric space. The surface exhibits visible jump discontinuities, which are induced by the sample-dependent mesh selection.}}
\end{subfigure}
\\
\begin{subfigure}{\linewidth}
	\centering
	\includegraphics[width=0.78\textwidth]{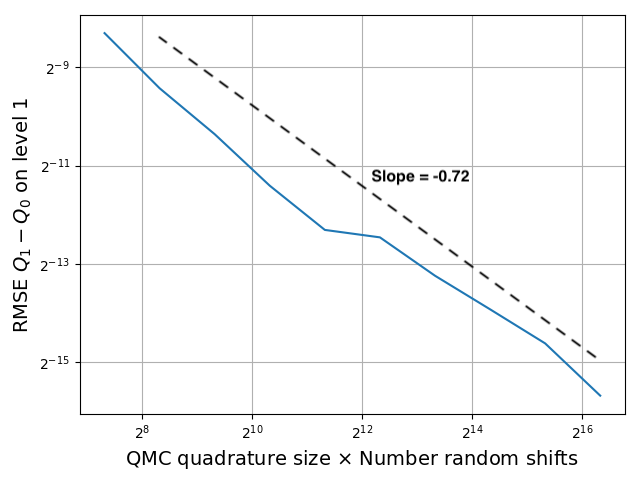}
	\caption{\revRtwo{The root-mean-square-error (RMSE) of QMC approximation of $\E{\abs*{{Q}(\mathbf{y}; \tol_1) - {Q}(\mathbf{y}; \tol_0)}}$. The convergence rate is reduced due to the discontinuities in the integrand. }}
\end{subfigure}
\caption{\revRtwo{The discontinuity induced by the sample-dependent mesh selection and its impact on the QMC convergence.}}
\label{fig:mesh_selection_discontinuity}
\end{figure}

The work~\cite{adMLMC_our} introduced the following sample-dependent quasi-optimal mesh selection strategy:
Let $u_h$ be a finite element solution of order $p$ to~\eqref{eq:bvp_general_1} on a mesh parametrized by $h$. 
The bias is assumed to be approximated by the leading order term in an error expansion, 
\begin{align}
	\label{eq:bias_model}
	Q(u(\mathbf{y})) - {Q}(u_h(\mathbf{y})) \approx \int_\mathcal{D} {\rho (\mathbf{x};\mathbf{y}) h^p(\mathbf{x};\mathbf{y})} d\mathbf{x},
\end{align}
where $\rho\in L_\mathbb{P}^{\frac{d}{p+d}}(\mathcal{D}\times\Omega)$ represents the error density, 
and $h: \mathcal{D} \times \Omega \to \mathbb{R}$ denotes the mesh-size function (see Theorem~2.1 in~\cite{adFEM_our}). 
This approximation is introduced to construct error estimates. As discussed in~\cite{adMLMC_our}, given a bias tolerance $\textrm{TOL}_{\textrm{bias}}$, with the optimal mesh function $h^* = h^*(\cdot, \mathbf{y})$, for each fixed $\mathbf{y}$ the error estimate satisfies the following:
\begin{align}
\label{eq:sample_dependent_stopping}
\int_{\mathcal{D}} \rho(\mathbf{x}; \mathbf{y}) h^{\ast} (\mathbf{x};\mathbf{y})^p d\mathbf{x} = \mathrm{TOL}_{\textrm{bias}} \frac{\int_{\mathcal{D}} \rho(\mathbf{x};\mathbf{y})^{\frac{d}{p+d}}d\mathbf{x}}{\int_{\mathcal{D}} \mathbb{E} \left[ \rho^{\frac{d}{p+d}} \right] }.
\end{align}
To avoid the computational expense of solving for $h^*(\mathbf{y})$ for each $\mathbf{y}$, we approximate the criterion~\eqref{eq:sample_dependent_stopping}
by restricting ourselves to a sequence of predesigned $h$-adaptive deterministic meshes \revRtwo{(generated by~\cite[Algorithm 2]{adMLMC_our}, illustrated by~\cite[Figure 2]{adMLMC_our}, and also discussed in Section~\ref{sec:numex_workflow})}, with mesh functions denoted by $h_k$, for $k \in \mathbb{N}_0$. 
Given a desired bias tolerance $\tol$ and a random sample $\mathbf{y} $, 
the mesh index $\mathcal{K}(\mathbf{y}) \in \mathbb{N}_0$ to evaluate $Q(u(\mathbf{y}))$ was chosen such that
\begin{align}
\label{eq:sample_dependent_stopping_discrete_deterministic}
\mathcal{K}(\mathbf{y}, \tol) =  \min \left\{ {k\in \mathbb{N}_0} ~\Bigg\vert~ \sum_K \rho_{k, K} (\mathbf y) h_{k, K}^{p+d} \leq \mathrm{TOL} \frac{\int_{\mathcal{D}} \rho(\mathbf x;\mathbf y)^{\frac{d}{p+d}}d \mathbf x}{\int_{\mathcal{D}} \mathbb{E} \left[ \rho^{\frac{d}{p+d}} \right] } \right\}.
\end{align}
\revAE{We emphasize the offline/online separation inherent in~\eqref{eq:sample_dependent_stopping_discrete_deterministic}: the predesigned mesh hierarchy $\{h_k\}_{k\in\mathbb{N}_0}$ and the normalizer $\int_{\mathcal D}\mathbb E[\rho^{d/(p+d)}]$ are constructed once offline (via~\cite[Algorithm~2]{adMLMC_our}) and reused across all parameter samples, whereas the per-sample indicators $\sum_K \rho_{k,K}(\mathbf y)\,h_{k,K}^{p+d}$ and the resulting stopping index $\mathcal{K}(\mathbf{y},\tol)$ are evaluated online for every realization of~$\mathbf y$. This split keeps the online cost dominated by FEM solves at the selected level and makes the online decision rule the only source of the sample-dependent discontinuities discussed next.}
This work extends the sampling method in~\cite{adMLMC_our} from MC to QMC. \revAE{However, a direct substitution is infeasible, 
because the sample-dependent mesh selection in~\eqref{eq:sample_dependent_stopping} introduces discontinuities in the 
parametric space and results in reduced convergence rates as demonstrated by Figure~\ref{fig:mesh_selection_discontinuity}.}

We therefore construct standard multilevel QMC estimators on these fixed adaptive meshes and analyze their computational complexity in the next section. Some other sophisticated methods addressing QMC regularity and adaptivity, including piecewise integration and the summation-by-parts technique, are discussed in~\ref{sec:adaptivity_qmc_piecewise_integration}--\ref{sec:summation_by_parts}.

\section{Monte Carlo and Quasi-Monte Carlo Estimators}
\label{sec:MC_QMC_estimators}
In this section, we introduce the MLMC and MLQMC estimators and discuss their computational complexities. \revAE{Table~\ref{tab:notation} collects the recurring symbols used in Sections~\ref{sec:prob}--\ref{sec:numex} for the reader's convenience.}

\revAE{%
\begin{table}[htbp]
\centering
\caption{\revAE{Notation used throughout Sections~\ref{sec:prob}--\ref{sec:numex}.}}
\label{tab:notation}
\begin{tabular}{ll}
\hline
Symbol & Meaning \\
\hline
$s$ & truncation dimension of the lognormal KL expansion \\
$\mathbf y = (y_1,\dots,y_s)$ & standard-Gaussian parameter vector \\
$\mathbf y_{-1} = (y_2,\dots,y_s)$ & parameter vector with the sign coordinate $y_1$ removed \\
$a(\mathbf x;\mathbf y)$ & lognormal diffusivity, Eq.~\eqref{eq:a_def} \\
$u(\mathbf x;\mathbf y)$ & pathwise solution of~\eqref{eq:bvp_general_1} \\
$Q,\ \tilde Q$ & QoI and its $y_1$-marginalized form, Section~\ref{sec:prob} \\
$h_k(\mathbf x)$ & deterministic adaptive mesh-size function on level $k$ \\
$\mathcal{K}(\mathbf{y},\tol)$ & sample-dependent stopping index, Eq.~\eqref{eq:sample_dependent_stopping_discrete_deterministic} \\
$\tol,\ \tol_\ell$ & global tolerance and per-level tolerance \\
$\alpha,\beta,\gamma$ & bias / variance / per-dimension cost exponents (Assumption~\ref{assumption:mc_complexity}) \\
$\workExp \coloneqq \gamma d$ & effective work-growth exponent, $C_k = \mathcal{O}(2^{\workExp k})$ \\
$\lambda$ & RQMC variance-decay exponent, Eq.~\eqref{eq:rqmc_variance_convergence} \\
$N_k,\ R$ & per-level QMC sample size and number of randomizations \\
$K_0,\ K_0^{\mathrm{Q}}$ & level-0 control-variate reduction factors (Section~\ref{sec:variance_reduction}) \\
\hline
\end{tabular}
\end{table}%
}

\subsection{Multilevel Monte Carlo Estimator}
We consider the MLMC estimator based on the $h$-adaptive meshes, where mesh functions are denoted by $h_k$, $k \in \mathbb{N}$. Let $\tilde{Q}_k$ be the QoI computed on the $k$-th mesh using the FEM method, and define
\begin{equation*}
	\Delta \tilde{Q}_k = \begin{cases}
	\tilde{Q}_k - \tilde{Q}_{k-1}, & \text{ $k \geq 1$},\\
	\tilde{Q}_0, & \text{ $k = 0$}.
	\end{cases}
\end{equation*} 
The MLMC estimator of $\mathbb{E}[\tilde{Q}]$ is then written as:
\begin{align}
\label{eq:biased_estimator_multilevel_type_i}
\hat{Q}_{\textrm{MLMC}} = \sum_{k = 0}^{K} \frac{1}{N_k } \sum_{i=1}^{N_k} \revAE{\Delta \tilde{Q}_k} (\mathbf{y}_{-1} (\omega_{i,k}) ),
\end{align}
where $N_k$ is the number of samples on level $k$, and $\mathbf{y}_{-1}(\omega_{i,k})$ are the independent and identically distributed (i.i.d.) random samples of $\mathbf{y}_{-1}$ across indices $i$ and levels $k$. The total number of levels $K$ is determined so that the bias satisfies:
\begin{equation}
	\abs*{\mathbb{E}[{\hat{Q}_{\textrm{MLMC}} - \tilde{Q}}]} \leq \tol_{\textrm{bias}}.
\end{equation}

Next, we list the standard MLMC assumptions.
\begin{assumption}[MLMC Assumptions]
	\label{assumption:mc_complexity}
	With the cost of computing one realization of $\tilde{Q}_k$ denoted by $C_k$, 
	assume that
	\begin{align}
		\abs*{\E{  {\tilde{Q}_k} - \tilde{Q}}} & =\mathcal{O} (2^{-\alpha k}),\\
		\E{\left( \Delta  \tilde{Q}_k \right)^2} & =\mathcal{O} (2^{-\beta k}),\\
		C_k & =\mathcal{O} \left(2^{k\workExp}\right),
	\end{align}
	as $k \to \infty$, where $\alpha \geq \frac{1}{2} \min(\beta, \workExp)$. {Here $\alpha$, $\beta$, and \revAE{$\workExp$} are the exponents governing bias decay, variance decay and cost growth, respectively. }
\end{assumption}
\revAE{To avoid the symbol clash between the per-dimension cost exponent~$\gamma$ and its combined form, we set $\workExp \coloneqq \gamma d$ once and for all, and refer to $\workExp$ as the \emph{effective work-growth exponent}. The values of $\workExp$ are regime-dependent: in the 2-D slit benchmark of Section~\ref{sec:numex}, uniform meshes give $\workExp = 2$ and goal-oriented $h$-adaptive meshes give $\workExp = 1$.}

Under Assumption~\ref{assumption:mc_complexity}, and by splitting the mean squared error equally between bias and statistical error, \revAE{i.e.,
\begin{align}
	\label{eq:mlmc_bias_error}
\left( \E{ {\hat{Q}_{\textrm{MLMC}} - \tilde{Q}}} \right)^2  &\leq \frac{1}{2} \tol^2,\\
\label{eq:mlmc_statistical_error}
\textrm{Var} \left[{\hat{Q}_{\textrm{MLMC}}}\right] &\leq \frac{1}{2} \tol^2,
\end{align}}
the overall MLMC work satisfies the following as $\tol \to 0$:
\begin{align}
\label{eq:complexity_SLMC}
\work_{\textrm{MLMC}} & =
\begin{cases}
\mathcal{O}(\tol^{-2}), & \text{if $\beta>\workExp$},\\
\mathcal{O}(\tol^{-2}\left(\log{\tol^{-1}}\right)^2), & \text{if $\beta=\workExp$},\\
\mathcal{O}(\tol^{-2\left(1+\frac{\workExp-\beta}{2\alpha}\right)}), & \text{if $\beta<\workExp$},
\end{cases}
&&  \text{as $\tol \to 0$}.
\end{align}
Details can be found in~\cite{AcNum_MLMC}.
For the numerical example with geometry-induced singularity that we consider in this work, using uniform meshes in $d=2$, we have $\beta = 2$, $\gamma = 1$. 
From~\eqref{eq:complexity_SLMC}, we obtain
\begin{align}
	W_{\textrm{MLMC}} = \mathcal{O}\!\left(\tol^{-2} \revAE{\left(\log{\tol^{-1}}\right)^2} \right), && \text{as $\tol \to 0$}.
\end{align}

\subsection{Multilevel Quasi-Monte Carlo Estimator}
First, we introduce notation for the QMC methods. The QMC estimator for the QoI, $\mathbb{E}[\tilde{Q}]$, is given by
\begin{align}
\hat{I}_N(\tilde{Q}) \coloneqq \frac{1}{N} \sum_{i=1}^N \tilde{Q}\left(\mathbf{y}(\mathbf{t}_i) \right).
\label{eq:qmc_estimation}
\end{align}
In integration w.r.t. the standard Gaussian measure, one maps each component of a low-discrepancy point $(\mathbf{t}_i)_j \in [0, 1]$ into $\mathbb{R}$ by
\begin{equation*}
	(\mathbf{y}(\mathbf{t}_i))_j \coloneqq \Phi^{-1}((\mathbf{t}_i)_j),\quad j = 1, \dotsc, s,
\end{equation*}
where $\Phi^{-1}$ is the inverse cumulative distribution function (CDF) of the standard Gaussian distribution, and $\{\mathbf{t}_i\}$, $i = 1, \dotsc, N$ is a predesigned deterministic low-discrepancy sequence in $[0, 1]^s$ (see~\cite{dick2010digital,niederreiter1992random}). However, using a deterministic point set introduces bias. Randomization techniques are introduced to address this, leading to the RQMC unbiased estimator:
\begin{equation}
	\hat{I}_N(\tilde{Q};\bm{\Delta}_r) = \frac{1}{N} \sum_{i=1}^N \tilde{Q}(\mathbf{y} (\mathbf{t}_i \oplus \bm{\Delta}_r)),
\end{equation}
where $\mathbf{t}_i$ denotes the $i$th deterministic QMC quadrature point, $\bm{\Delta}_r$ represents the $r$-th randomization, and $\oplus$ denotes the randomization operation. An example of such an operation is the random shift, where $\bm{\Delta}_r \sim U[0, 1]^s$ and $a \oplus b = (a + b) \textrm{ mod } 1$, with the modulo taken componentwise. 
Typically, $R \ll N$ randomizations are used for a practical variance estimate as follows:
\begin{align}
\hat{I}_{N,R} (\tilde{Q}) \coloneqq \frac{1}{R} \sum_{r=1}^{R} \frac{1}{N} \sum_{i=1}^N \tilde{Q}(\mathbf{y} (\mathbf{t}_i \oplus \bm{\Delta}_r)).
\end{align}

\subsection{Quasi-Monte Carlo Adaptive Finite Element Estimator}
Similarly to the MLMC estimator $\hat{Q}_{\textrm{MLMC}}$, we define the MLQMC estimator by:
\begin{align}
	\begin{split}
	\hat{Q}_{\textrm{MLQMC}} &= \sum_{k = 0}^{K} \frac{1}{R}\sum_{r=0}^{R-1} \frac{1}{N_k } \sum_{j=0}^{N_k-1}  \Delta \tilde{Q}_{k} (\mathbf{y}_{-1} (\mathbf{t}_j^{k,r}))\\
	&= \sum_{k=0}^{K} \frac{1}{R} \sum_{r=0}^{R - 1} I_{N_k}(\Delta \tilde{Q}_{k}; \bm{\Delta}^k_r),
	\end{split}
\end{align}
with
\begin{equation}
	I_{N_k}(\Delta \tilde{Q}_{k}; \bm{\Delta}^k_r) \coloneqq \frac{1}{N_k} \sum_{j=0}^{N_k - 1} {\Delta \tilde{Q}_k(\mathbf{y}_{-{1}}(\mathbf{t}_j \oplus \bm{\Delta}_{r}^{k} ))}.
\end{equation}
The randomizations $\bm{\Delta}^k_r$ are drawn from a uniform distribution on $[0, 1]^s$, i.i.d. for each level $k$ and randomization $r$. \revAE{Throughout this work we instantiate the point set $\{\mathbf t_j^{k,r}\}$ using scrambled Sobol' sequences; $\bm{\Delta}^k_r$ then denotes the $r$-th random scramble on level~$k$. Other admissible RQMC constructions (such as randomly shifted rank-1 lattice rules, as described in~\cite{KSSSU2017}) would fit the same abstract framework.}

We now examine the variance convergence of the RQMC estimator level-wise. In standard QMC theory, for each level $k = 0, 1, \dotsc, K$, one assumes the existence of constants $V^{\mathrm{Q}}_{k} > 0$ and a convergence rate $\lambda > 1$ such that, for sufficiently large $N_k$,
\begin{equation}
\label{eq:rqmc_variance_convergence}
\textrm{Var} \left[{I_{N_k} (  \Delta \tilde{Q}_{k};\bm{\Delta} ) }\right] \leq \frac{V^{\mathrm{Q}}_{k}}{N_k^\lambda}.
\end{equation}
The work~\cite{liu2023nonasymptotic} shows that, for the RPDE model considered in this work, the asymptotic convergence rate is $\lambda = 2 - \epsilon$ (for any $\epsilon > 0$) with {$V^{\mathrm{Q} }_{k} \coloneqq V^{\mathrm{Q} }_{k}(\epsilon)$, and $\lim_{\epsilon \to 0} V^{\mathrm{Q} }_{k}(\epsilon) = +\infty$.}

However, \revAE{in the finite-sample regime our numerical results exhibit a smaller effective $\lambda < 2$. On the 2-D slit benchmark of Section~\ref{sec:numex} we measure an observed slope $\lambda \approx 1.6$ for the level-wise RQMC variance $\mathrm{Var}[I_{N_k}(\Delta\tilde Q_k;\bm\Delta)]$ as $N_k$ is increased; see Figure~\ref{fig:ex_2_qmc_convergence_is} (left panel, without importance sampling) and the accompanying discussion in Section~\ref{sec:numex}. This pre-asymptotic rate degradation from the theoretical $\lambda = 2-\epsilon$ is consistent with~\cite{liu2023nonasymptotic}. }

In the RQMC method, one uses $N_k$ samples and $R$ randomizations for each level $k$ to obtain a practical variance estimate, although a single randomization suffices to provide an unbiased result. The total computational cost is given by
\begin{align}
	R \sum_{k=0}^{K} C_k N_k.
\end{align}
Following the assumption~\eqref{eq:rqmc_variance_convergence}, the variance of the RQMC estimator can be bounded by:
\begin{align}
	\sum_{k=0}^{K} \frac{1}{R} \textrm{Var}  \left[{I_{N_k} (\Delta \tilde{Q}_{k}; \bm{\Delta}) }\right] \leq \frac{1}{R} \sum_{k = 0}^{K} \frac{V_k^{\mathrm{Q}}}{N_k^\lambda}.
\end{align}
Following~\cite{KSSSU2017}, we list the standard MLQMC assumptions:
\begin{assumption}[MLQMC Assumptions]
	\label{assumption:qmc_rates}
	There exist constants $\alpha, \beta, \revAE{\workExp} > 0$ and $\lambda > 1$, such that as $k \to \infty$,
	\begin{align}
	\begin{split}
		\abs*{\E{ \tilde{Q}_k - \tilde{Q}}} &= \mathcal{O} \left( 2^{-\alpha k} \right)\\
		\mathrm{Var}  \left[{I_{N_k} ( \Delta \tilde{Q}_{k}; \bm{\Delta}_r^{k}) }\right] &= \mathcal{O} \left(2^{-\beta k} N_k^{-\lambda} \right)\\
		C_k &= \mathcal{O} \left(2^{k\workExp}\right),
	\end{split}
	\end{align}
  \revAE{with $\alpha \geq \frac{1}{2} \min(\beta, \workExp)$}.
	Here, $\alpha$, $\beta$, and \revRtwo{$\workExp=\gamma d$ have} 
	the same meaning as in Assumption~\ref{assumption:mc_complexity}, 
	and $\lambda = 2 - \epsilon$, for any $\epsilon > 0$ to be specified later, is the RQMC variance decay exponent from~\eqref{eq:rqmc_variance_convergence}.
\end{assumption}
\revRtwo{In the 2-D slit problem of Section~\ref{sec:numex}, with the particular particular sequences of refinements used there, uniform meshes give $\workExp = 2$ and goal-oriented $h$-adaptive meshes give $\workExp = 1$, while uniform and adaptive meshes give the same $\alpha$ and $\beta$, see the discussion following~\eqref{eq:adaptive_work_scaling}.}

If Assumption~\ref{assumption:qmc_rates} holds and the bias and statistical errors are split equally, the MLQMC complexity satisfies:
\begin{align}
\label{eq:complexity_MLMC}
\work_{\textrm{MLQMC}} &=
\begin{cases}
\mathcal{O} \left(\tol^{-2/\lambda}\right), & \text{if $\beta> \lambda\,\workExp$},\\
\mathcal{O} \left(\tol^{-2/\lambda}\left(\log{\tol^{-1}}\right)^{1/\lambda + 1}\right), & \text{if $\beta=\lambda\,\workExp$},\\
\mathcal{O} \left(\tol^{-2/\lambda - \left(\frac{\lambda\,\workExp -\beta}{\alpha \lambda}\right)}\right), & \text{if $\beta<\lambda\,\workExp$}.
\end{cases}
\end{align}
Detailed derivations of the complexity can be found in~\cite{KSSSU2017}. 

Notice that the variance $\mathrm{Var}  \left[{I_{N_k} ( \Delta \tilde{Q}_{k}; \bm{\Delta}_r^{k}) }\right]$ decays with both the QMC quadrature size $N_k$ and the level parameter $k$. For our adaptive meshes, we have
\begin{equation}
	(Q_{k} - Q_{k - 1})(\mathbf{y}) \simeq \int_{\mathcal{D}} \rho(\mathbf{x};\mathbf{y}) \left( h_{k}^p(\mathbf{x}) - h_{k - 1}^p(\mathbf{x}) \right) d\mathbf{x},
\end{equation}
where $h_{k}(\mathbf{x})$, $k = 0, 1, \dotsc$, is the deterministic mesh function on level $k$, determined by
\begin{equation*}
	h_{k} (\mathbf x) = \tol_{k}^{1/p} \frac{\rho_0^{-\frac{1}{p+d}}(\mathbf x)}{\left( \int_{\mathcal{D}} \rho_0^{\frac{d}{p+d}}(\mathbf x) d\mathbf x \right)^{1/p}},
\end{equation*}
where $\rho_0$ is the error density function with a constant diffusion coefficient $a \equiv 1$. \revRtwo{The specific choice $\tol_k = 2^{-k-2}$ as in~\cite{adMLMC_our} is used in the numerical experiments of Section~\ref{sec:numex} and documented in Section~\ref{sec:numex_workflow}.}
{
\revRtwo{Let $C := \tol_{k - 1}/\tol_{k}$ denote the inter-level tolerance ratio; with the choice above, $C = 2$ is a fixed constant independent of~$k$.} Thus, we have
\begin{equation}
	(Q_{k} - Q_{k - 1})(\mathbf{y}) \simeq \tol_{k}(C^{-1} - 1) \int_{\mathcal{D}} \rho(\mathbf{x};\mathbf{y}) \frac{\rho_0^{-\frac{p}{p+d}}(\mathbf x)}{\left( \int_{\mathcal{D}} \rho_0^{\frac{d}{p+d}}(\mathbf x) d\mathbf x \right)} d\mathbf{x}.
\end{equation}

The mean and variance of $(Q_{k} - Q_{k - 1})$ exhibit the desired decay rates of $\mathcal{O}(\tol_{k})$ and $\mathcal{O}(\tol_{k}^2)$, respectively, provided that $\E{\left(\int_{\mathcal{D}} \rho(\mathbf{x};\mathbf{y}) \rho_0^{-\frac{p}{p+d}}(\mathbf x) d\mathbf{x} \right)^2} < +\infty$. This decay behavior is corroborated by the numerical experiments presented in Section~\ref{sec:numex}. 
}

The decay of $\mathrm{Var}  \left[{I_{N_k} ( \Delta \tilde{Q}_{k}; \bm{\Delta}_r^{k}) }\right]$ with respect to the QMC quadrature size $N_k$ can be analyzed by examining the derivatives of $\Delta \tilde{Q}_{k}$ with respect to $\mathbf{y}$~\cite{liu2023nonasymptotic}. In the following we show the derivatives of $Q - Q_h$ for notation simplicity, and similar results extend to $\Delta \tilde{Q}_{k}$. Following~\cite{liu2023nonasymptotic}, we have
\begin{align}
		\abs{\partial^{\mathfrak{u}} Q(y)} &\leq \norm{Q}_{V^{\prime}} \norm{f}_{V^{\prime}} \frac{\abs{\mathfrak{u}}!}{(\log 2)^{\abs{\mathfrak{u}}}} \left( \prod_{j \in \mathfrak{u}} b_j \right) \left( \prod_{j=1}^s \exp (b_j \abs{y_j}) \right)\\
		\abs{\partial^{\mathfrak{u}} Q_h(y)} &\leq \norm{Q}_{V^{\prime}} \norm{f}_{V_h^{\prime}} \frac{\abs{\mathfrak{u}}!}{(\log 2)^{\abs{\mathfrak{u}}}} \left( \prod_{j \in \mathfrak{u}} b_j \right) \left( \prod_{j=1}^s \exp (b_j \abs{y_j}) \right).
\end{align}
This yields:
\begin{equation}
	\abs{\partial^{\mathfrak{u}} (Q(y) - Q_h(y))} \leq \left( \norm{Q}_{V^{\prime}} \norm{f}_{V^{\prime}} + \norm{Q_h}_{V_h^{\prime}} \norm{f}_{V_h^{\prime}} \right) \frac{\abs{\mathfrak{u}}!}{(\log 2)^{\abs{\mathfrak{u}}}} \left( \prod_{j \in \mathfrak{u}} b_j \right) \left( \prod_{j=1}^s \exp (b_j \abs{y_j}) \right).
\end{equation}
We employ a transformation-based approach to map a QMC quadrature $\mathbf{t} \in [0, 1]^s$ to $\mathbf y \in \mathbb{R}^s$ via $\mathbf{y} = \Phi^{-1}(\mathbf{t})$, where $\Phi^{-1}$ denotes the inverse CDF of the standard Gaussian distribution. We have the equivalence of the integrals:
\begin{equation}
	\int_{\mathbb{R}^s} Q(y) \varphi(\mathbf{y}) d\mathbf{y} = \int_{[0, 1]^s} Q(\Phi^{-1}(\mathbf{t})) d\mathbf{t}. 	
\end{equation}
The derivatives of $Q$ with respect to $\mathbf{t}$ are then given by:
\begin{equation}
	\revAE{\abs*{\frac{\partial^{\mathfrak{u} } Q}{\partial \bm{t}_{\mathfrak{u}}}}}  = \abs{\partial^{\mathfrak{u}} Q(y)} \prod_{j\in \mathfrak{u}} \abs{\partial^{j} \Phi^{-1}(t_j)}.
\end{equation}
Following the derivations in~\cite{liu2023nonasymptotic}, we have
\begin{equation}
	\abs{\partial^{j} \Phi^{-1}(t_j)} \leq C \min(t_j, 1 - t_j)^{-1},
\end{equation}
for a given constant $C > 0$ and $j = 1, \dotsc, s$. Thus,
\begin{equation}
	\label{eq:derivative_Q_minus_Qh}
	\abs*{\frac{\partial^{\mathfrak{u}} \left(Q - Q_h\right) }{\partial \mathbf{t}_{\mathfrak{u}}} } \leq C_{s, \delta} \prod_{j=1}^s \min(t_j, 1 - t_j)^{-\mathbb{I}_{j \in \mathfrak{u}} - \delta},  
\end{equation}
where $\delta > 0$ can be chosen arbitrarily small, and $C_{s, \delta}$ is a constant depending on the dimension $s$ and $\delta$. Equation~\eqref{eq:derivative_Q_minus_Qh} implies a variance decay rate $\lambda = 2 - \epsilon$ for any $\epsilon > 0$ in the RQMC variance decay assumption~\eqref{eq:rqmc_variance_convergence}, following the derivations in~\cite{liu2023nonasymptotic}.

The work model for adaptive meshes is given by:
\begin{equation}
	\label{eq:adaptive_work_scaling}
	\revRtwo{C_{k} = \int_{\mathcal{D}} h_{k}^{-d}(\mathbf x)\, d\mathbf{x} = \mathcal{O}\!\left(\tol_{k}^{-d/p}\right),}
\end{equation}
\revRtwo{where the suppressed prefactor depends only on $\rho_0$ and the domain $\mathcal{D}$; the full algebraic expression is reported in~\cite{adMLMC_our}. Writing this  cost per level in the standard MLMC parametrization~\eqref{eq:complexity_MLMC}, $C_k = \mathcal{O}(2^{k\workExp})$, one obtains an effective work-growth exponent $\workExp = \gamma d$ that depends on whether the mesh hierarchy is uniform or adaptive: uniform meshes on the 2-D slit problem yield $\gamma = 1$ and therefore $\workExp = 2$, whereas the goal-oriented adaptive hierarchy of~\cite{adFEM_our,adMLMC_our} concentrates elements near the re-entrant corner and delivers an effective $\gamma = 1/2$ (hence $\workExp = 1$); see Figure~\ref{fig:single_level_k_mc}(c) for the empirical slope of $C_k$ versus~$k$. The adaptive value of~$\workExp$ is therefore \emph{a property of the singular geometry together with the goal-oriented adaptive refinement strategy}, not of the finite element approximation order~$p$ per se.}

In the following we compare the complexities of the MLQMC method for uniform and adaptive meshes on the 2-D slit problem as studied in~\cite{adMLMC_our}, to be specified in Section~\ref{sec:numex}.
For uniform meshes, \revAE{due to the geometric singularity,} we have $\alpha = 1$, $\beta = 2$, and $\gamma  = 1$ (so that $\workExp = \gamma d = 2$). With $\lambda = 2 - \epsilon$ for any arbitrarily small and fixed $\epsilon > 0$, we find $\beta < \lambda\,\workExp$, yielding the complexity:
\begin{align}
	W_{\textrm{MLQMC}} = \mathcal{O} \left( \tol^{-2} \right).
\end{align}

In contrast, for adaptive meshes in 2-D slit problem, we have that \revRtwo{$\alpha = 1$, $\beta = 2$, $\gamma  = 1/2$ (so that $\workExp = \gamma d = 1$)}, which are also observed in Figure~\ref{fig:single_level_k_mc} of Section~\ref{sec:numex}.

With $\lambda = 2 - \epsilon$ {for any fixed} and arbitrarily small $\epsilon > 0$, we find that $\beta > \lambda\,\workExp$, yielding,
\begin{align}
\label{eq:mlqmc_adaptive_2d_complexity}
\revRtwo{W_{\textrm{MLQMC}} = \mathcal{O}\!\left( \tol^{-2/(2-\epsilon)} \right), \qquad \epsilon > 0,}
\end{align}
\revAE{which is arbitrarily close to $\mathcal{O}(\tol^{-1})$ and, crucially, strictly better than the $\mathcal{O}(\tol^{-2})$ rate attained by the uniform-mesh MLQMC estimator. Notice that one cannot take $\epsilon \to 0$ since the implicit constant diverges to $+\infty$, as discussed after~\eqref{eq:rqmc_variance_convergence}.} The complexities of the MLQMC method for higher physical dimensions $d \geq 3$ can be derived similarly and are summarized in Table~\ref{tab:complexity_summary}.

\begin{table}[htbp]
\centering
\caption{\revRtwo{The complexities of MLQMC of the slit problem for uniform and adaptive meshes, with $\lambda = 2-\epsilon$, arbitrarily small $\epsilon > 0$. With uniform meshes, all dimension settings ($d \ge 2$) correspond to the worst regime $\beta < \lambda\,\workExp$. The 2-D adaptive row lies in the best regime $\beta > \lambda\,\workExp$ of Eq.~\eqref{eq:complexity_MLMC}; the remaining rows correspond to the worst regime $\beta < \lambda\,\workExp$. }}
\begin{tabular}{ccc}
\hline
Dimension &  Uniform Mesh &  Adaptive Mesh\\
\hline
2-D MLQMC & $\mathcal{O}\left(\tol^{-2}\right)$ & \revRtwo{$\mathcal{O}\left(\tol^{-2/(2-\epsilon)}\right), \epsilon > 0$} \\
3-D MLQMC & $\mathcal{O}\left(\tol^{-3}\right)$ & $\mathcal{O}\left(\tol^{-1.5}\right)$\\
$\vdots$ & & \\
d-D MLQMC & $\mathcal{O}\left(\tol^{-d}\right)$ & $\mathcal{O}\left(\tol^{-d/2}\right)$\\
\hline
\end{tabular}
\label{tab:complexity_summary}
\end{table}

In this section we have discussed the complexity of the MLQMC algorithm under the assumptions of the convergence rates of the bias and variance. To further enhance the performance of the MLQMC algorithm, we explore variance reduction techniques, i.e., the importance sampling and control variate within the MLQMC framework in the next section. 
\section{Two approaches for variance reduction}
\label{sec:variance_reduction}

This section introduces two variance reduction strategies within the multilevel hierarchy. The first is importance sampling (IS), which reduces level-wise variance by changing the integration measure; the second is the control variate, which targets variance reduction on level 0, the level that typically dominates the overall computational cost across all levels.
\subsection{Importance Sampling}
\label{sec:theory_is}

Following the work~\cite{liu2023nonasymptotic}, we introduce a Gaussian proposal $\varphi_{\bm{\alpha}}$, 
\begin{equation}
\varphi_{\bm{\alpha}}(\mathbf{y}) = \prod_{j=1}^s \varphi_{\alpha_j}(y_j) =  \prod_{j=1}^s \frac{1}{{\alpha}_j \sqrt{2\pi} } \exp\left( -\frac{{y}_j^2}{2{\alpha}_j^2} \right),
\end{equation}
parametrized by $\bm{\alpha} = (\alpha_1, \dotsc, \alpha_s)$. This proposal targets integration w.r.t. the standard Gaussian distribution. The singularity at the boundary can be mitigated when ${\alpha}_j > 1$ for all $ j = 1, \dotsc, s$. 

The expectation of the QoI can be written as
\begin{align}
	\label{eq:is_expectation}
	\begin{split}
		\mathbb{E}_{\varphi}[Q] &= \int_{\mathbb{R}^s} Q(\mathbf{y}) \varphi(\mathbf{y}) d\mathbf{y}\\
		&= \int_{\mathbb{R}^s} Q(\mathbf{y}) \frac{\varphi( \mathbf{y})}{\varphi_{\boldsymbol{\alpha} }(\mathbf{y})} \varphi_{\boldsymbol{\alpha}}(\mathbf{y}) d\mathbf{y}\\
		&= \mathbb{E}_{\varphi_{\boldsymbol{\alpha} }}\left[ Q_{\textrm{IS}}  \right],
	\end{split}
\end{align}
where we denote $Q_{\textrm{IS}}(\mathbf{y}) \coloneqq Q(\mathbf{y}) \frac{\varphi( \mathbf{y})}{\varphi_{\boldsymbol{\alpha} }(\mathbf{y})}$. Notice that~\eqref{eq:is_expectation} can also be equivalently formulated by a change of variables. To see this, we denote $\mathbf{y} = \Phi^{-1}_{\bm{\alpha}}(\mathbf{t})$, where ${y}_j = \Phi^{-1}_{{\alpha}_j}({t}_j)$ and $\Phi_{{\alpha}_j}$ is the CDF of the proposal distribution $\varphi_{{\alpha}_j}$\revAE{; consequently $\mathbf{t} = \Phi_{\bm{\alpha}}(\mathbf{y}) \in [0,1]^s$, so that the integration domain below is the unit cube rather than $\mathbb{R}^s$.} Then we have
\begin{equation}
	\label{eq:is_change_of_variables}
	\begin{split}
		\mathbb{E}_{\varphi}[Q] &= \revAE{\int_{[0,1]^s}} Q(\Phi^{-1}_{\bm{\alpha}}(\mathbf{t})) \varphi(\Phi^{-1}_{\bm{\alpha}}(\mathbf{t})) \left| \det \nabla \Phi^{-1}_{\bm{\alpha}}(\mathbf{t}) \right| d\mathbf{t}\\
		&= \revAE{\int_{[0,1]^s}} Q(\Phi^{-1}_{\bm{\alpha}}(\mathbf{t})) \frac{\varphi(\Phi^{-1}_{\bm{\alpha}}(\mathbf{t}))}{\varphi_{\bm{\alpha}} (\Phi^{-1}_{\bm{\alpha}}(\mathbf{t}) ) } d\mathbf{t}\\
		&= \mathbb{E}_{\varphi_{\boldsymbol{\alpha} }}\left[ Q_{\textrm{IS}}  \right],
	\end{split}
\end{equation}
where, in the second line, the inverse function theorem ensures the equivalence of the IS weight through the Jacobian determinant of $\Phi^{-1}_{\bm{\alpha}}$. \revAE{The rightmost expectation is thus an integral against the uniform measure on $[0,1]^s$ and is suitable for direct QMC evaluation.}

We select the quasi-optimal $\bm{\bar{\alpha}}$ as any minimizer of the  empirical second-order moment of $Q_{\textrm{IS}}$ based on $n$ samples:
\begin{align}
	\label{eq:bar_alpha_is_optimization}
\bar{\bm{\alpha}} \in \arg \min_{{\alpha}_j \geq 1, j = 1, \dotsc, s}  \frac{1}{n} \sum_{i=1}^n Q^2( \mathbf{y}_i) \cdot \frac{\varphi( \mathbf{y}_i)}{\varphi_{\boldsymbol{\alpha}} (\mathbf{y}_i)}.
\end{align}
\revAE{Here $\{\mathbf{y}_i\}_{i=1}^n$ is an i.i.d.\ training sample drawn from the \emph{target} standard Gaussian density~$\varphi$, and the constraint $\alpha_j \geq 1, j = 1, \dotsc, s$ ensures the integrand regularity is not compromised.}
In this work, we will consider the change of measure with respect to $\mathbf{y}_{-{1}}$, i.e. we fix ${\alpha}_1 = 1$ and consider ${\alpha}_j \geq 1$ for $j = 2, 3, \dotsc, s$. The connections of the importance sampling with the original optimization algorithm proposed in~\cite{adMLMC_our} will be discussed in the~\ref{sec:mlqmc_is_appendix}.

However, as we will observe in the numerical section, the importance sampling only yields moderate variance reduction of the MC/QMC estimator of $\tilde{Q}_{\ell} - \tilde{Q}_{\ell - 1}$ for the specific RPDE problem we have considered (see Figure~\ref{fig:ex_2_qmc_convergence_is}). Although the moderated measure improves the integrand regularity, it can inflate the integrand variance (see~\cite{liu2023nonasymptotic} for instance). Moreover, optimizing the proposal parameter $\bm{\alpha}$ becomes increasingly costly as $s$ increases. In the next section we will discuss an alternative approach, i.e., the control variate. 
\subsection{Control variate on level 0}
\label{sec:cv_level_0}
In scenarios where the cost is dominated by level 0 within the non-asymptotic regime of the multilevel hierarchy, the total work required for a given finite tolerance becomes indistinguishable between the uniform and adaptive approaches when they share the same coarsest mesh. An example of this behavior is observed in~\cite{adMLMC_our}. 

In this section, we introduce a control variate (CV) to mitigate the cost contribution at level 0. Notably, the study~\cite{Nobile2015multi} addresses a PDE problem with a rough random field and \revAE{implements} a control variate at each level of the multilevel hierarchy. In \cite{Nobile2015multi}, the decision to introduce a CV in each level is \revAE{contingent on} the smoothness of the field, denoted as $C^{\alpha}$\revAE{,} and it is beneficial when $\alpha \geq 1$. Specifically, the CV relies on a convolution-based smooth field, with its expectation computed via stochastic collocation to exploit the enhanced smoothness of the coefficient field.

However, this approach does not reduce the input dimension of the CV, making the collocation step expensive in high dimensions. Moreover, since our primary focus is on reducing costs at level 0, we limit the application of the CV to this initial level only. {We leave the extension of the CV to higher levels for future work.}

Below, we formalize the assumptions for the control variate at level 0.
\begin{assumption}[Level-0 CV effectiveness]
	\label{ass:cv_effective}
	We consider a control variate $\tilde{Q}^{\mathrm{CV}}_0$ with reduced input dimension compared to $\tilde{Q}_0$. This control variate is assumed to provide effective variance reduction for both MC and RQMC estimators, in the sense that
	\begin{equation}
		\var{\tilde{Q}_0 - \tilde{Q}_0^{\mathrm{CV}}} \leq K_0 \var{\tilde{Q}_0}
	\end{equation}
	and
	\begin{equation}
		\var{\tilde{Q}_0 - \tilde{Q}_0^{\mathrm{CV}}; \bm{\Delta}} \leq K_0^Q \var{\tilde{Q}_0;\bm{\Delta}}
	\end{equation}
	hold with constants $K_0,\ K_0^Q \ll 1$.
\end{assumption}

\begin{figure}[ht]
    \centering
    \begin{minipage}[t]{0.50\textwidth}
        \centering
        \includegraphics[width=\linewidth]{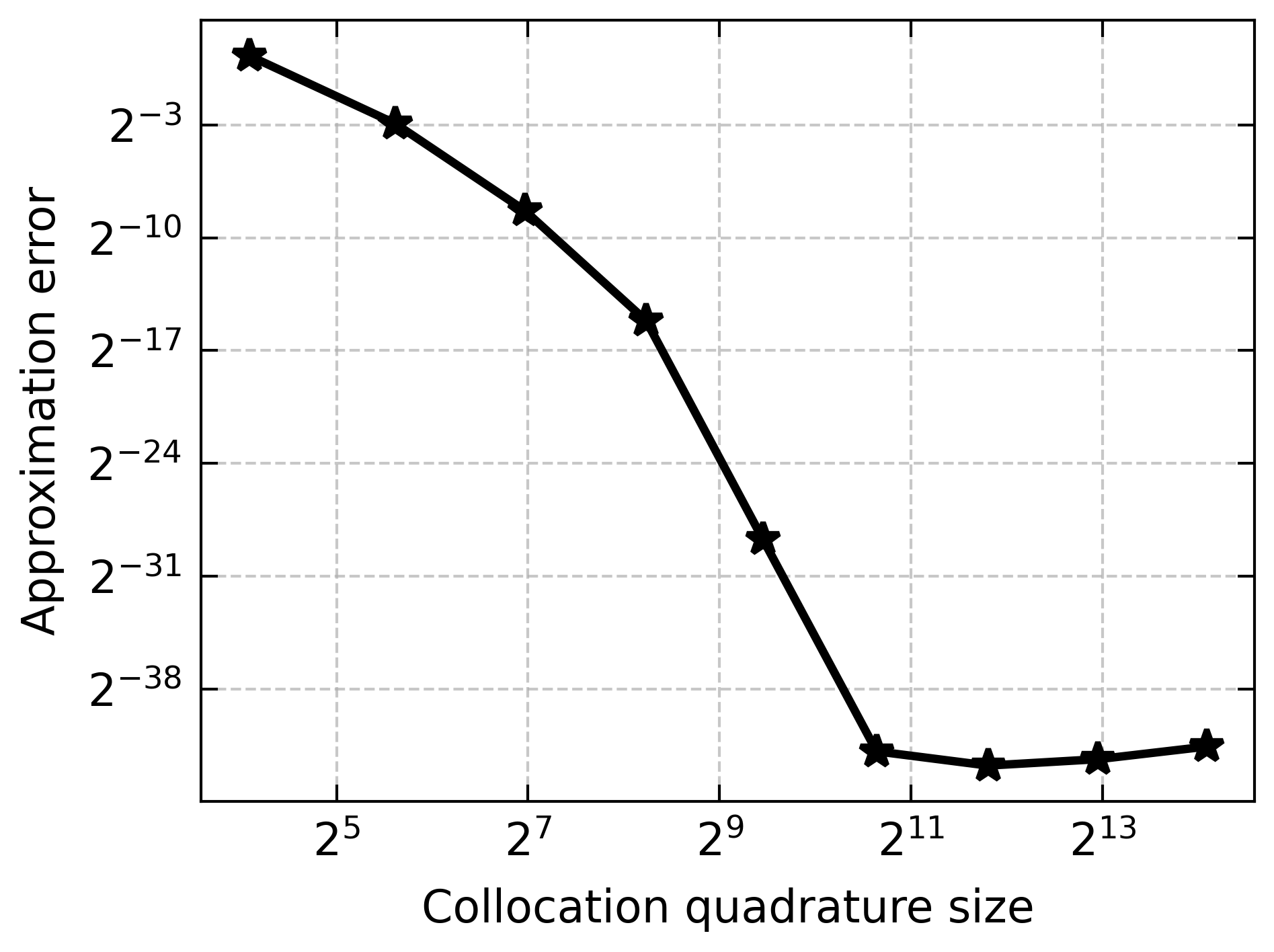} 
        \subcaption{The exponential decay of the stochastic collocation approximation error.}
        \label{fig:left}
    \end{minipage}
    \begin{minipage}[t]{0.48\textwidth}
		\vspace{-4.6cm}
        \centering
		\begin{subfigure}{0.48\linewidth}
			\includegraphics[width = \textwidth]{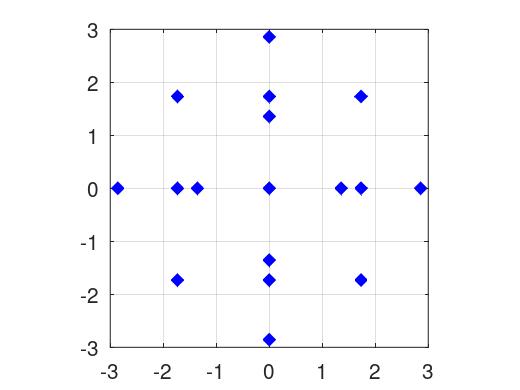}
			\caption{$w = 2$}\label{subfig:w=2}
		\end{subfigure}
		\begin{subfigure}{0.48\linewidth}
			\includegraphics[width = \textwidth]{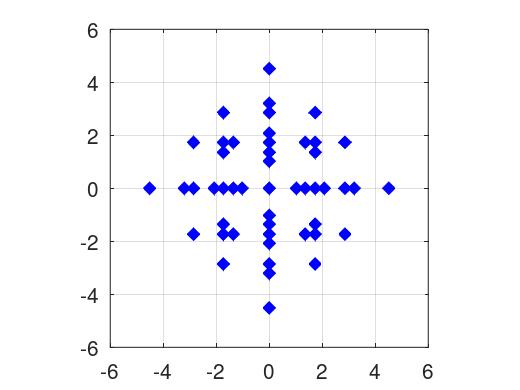}
			\caption{$w = 3$}\label{subfig:w=3}
		\end{subfigure}	
		\\
		\begin{subfigure}{0.48\linewidth}
			\includegraphics[width = \textwidth]{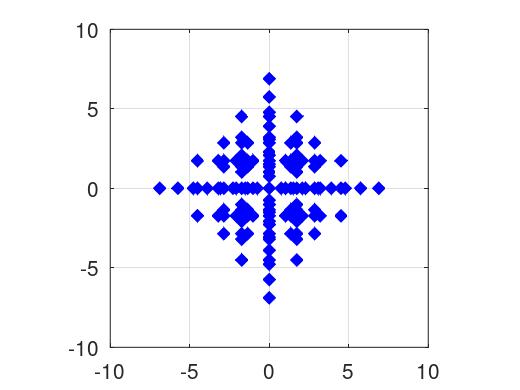}
			\caption{$w = 4$}\label{subfig:w=4}
		\end{subfigure}	
		\begin{subfigure}{0.48\linewidth}
			\includegraphics[width = \textwidth]{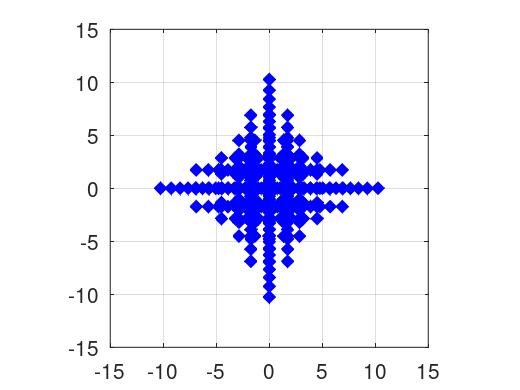}
			\caption{$w = 5$}\label{subfig:w=5}
		\end{subfigure}	
    \end{minipage}
	\hfill
    \caption{Stochastic collocation using Smolyak sparse grids with Gauss--Hermite quadrature in two dimensions. Quadrature levels $w = 2, 3, 4, 5$ are shown. The left figure displays approximation error convergence, while the right figures illustrate quadrature node distributions for each level $w$.}
    \label{fig:stochastic_collocation}
\end{figure}

\revAE{\paragraph{Numerical observation (stochastic collocation error)} Separately from the variance-reduction hypothesis of Assumption~\ref{ass:cv_effective}, we rely on the following empirical observation to justify using stochastic collocation for the CV expectation $\mathbb{E}[\tilde{Q}_0^{\mathrm{CV}}]$: on the 2-D slit benchmark of Section~\ref{sec:numex}, the stochastic collocation approximation of $\mathbb{E}[\tilde{Q}_0^{\mathrm{CV}}]$ converges exponentially in the Smolyak level~$w$ (Figure~\ref{fig:stochastic_collocation}, left panel), so the collocation error is negligible at the QMC tolerances targeted in Section~\ref{sec:numex}. We report this as a numerical observation rather than a hypothesis because it is specific to the rank-2 SVD CV and the 49-dimensional coefficient field of Section~\ref{sec:numex}, and we do not formally certify exponential decay for other CV families or random-field configurations.}



{Figure~\ref{fig:stochastic_collocation} displays the stochastic collocation approximation error of the CV based on a rank-2 SVD. The left plot shows an exponential decay, while the right plots illustrate the Gauss--Hermite quadrature nodes for levels $w = 2, 3, 4, 5$, with more nodes at higher levels.}

\begin{remark}[Optimal choice of level 0]
	The introduction of the CV could alter the optimal choice of level 0. As discussed in~\cite{AcNum_MLMC}, level 0, the coarsest model in the multilevel hierarchy, is beneficial to retain if
	\begin{equation}
		\sqrt{V_{0} C_{0}} + \sqrt{V_1 C_1} < \sqrt{\var{Q_1} C_1} \approx \sqrt{V_0 C_1},
		\label{eq:optimal_level_0_old}
	\end{equation}
	where \(V_0 = \var{Q_0}\) and \(V_1 = \var{Q_1 - Q_0}\). Otherwise, level 0 should be discarded and the current level 1 becomes the new ``level 0''. Similarly, for MLQMC, the optimal level 0 is retained if
	\begin{equation}
		\label{eq:optimal_level_0_old_qmc}
		{(V_{0})^{Q}}^{\frac{1}{\lambda + 1}} {C_{0}^{\frac{\lambda}{\lambda + 1}}} + {(V_1^{Q})}^{\frac{1}{\lambda + 1}} {C_1}^{\frac{\lambda}{\lambda + 1}} < {\var{Q_1;\bm{\Delta}} }^{\frac{1}{\lambda + 1}} {C_1}^{\frac{\lambda}{\lambda + 1}} \approx {(V_0^{Q})}^{\frac{1}{\lambda + 1}} {C_1}^{\frac{\lambda}{\lambda + 1}}.
	\end{equation}
	Under Assumption~\ref{ass:cv_effective} which incorporates the CV, we propose revised criteria to retain level 0 for MLMC and MLQMC as follows:
	\begin{subequations}
		\begin{equation}
			\sqrt{K_0 V_{0} C_{0}} + \sqrt{V_1 C_1} < \sqrt{K_1 V[Q_1] C_1} \approx \sqrt{K_1 V_0 C_1}.
			\label{eq:optimal_level_0_new_mc}
		\end{equation}
		\begin{equation}
			{(K_0^{Q} V_{0}^{Q})}^{\frac{1}{\lambda + 1}} {C_{0}}^{\frac{\lambda}{\lambda + 1}} + {(V_1)^{Q}}^{\frac{1}{\lambda + 1}} {C_1}^{\frac{\lambda}{\lambda + 1}} < (K_1^{Q} \var{Q_1;\bm{\Delta}})^{\frac{1}{\lambda + 1}} {C_1}^{\frac{\lambda}{\lambda + 1}} \approx (K_1^{Q} V_0^{Q})^{\frac{1}{\lambda + 1}} {C_1}^{\frac{\lambda}{\lambda + 1}},
			\label{eq:optimal_level_0_new_qmc}
		\end{equation}
		\label{eq:optimal_level_0_new}
	\end{subequations}
    where \(K_0, K_1, K_0^Q, K_1^Q \ll 1\) are variance reduction factors. Notably, the previously optimal level 0 satisfying~\eqref{eq:optimal_level_0_old} and~\eqref{eq:optimal_level_0_old_qmc}  might no longer satisfy the new criteria~\eqref{eq:optimal_level_0_new_mc} and~\eqref{eq:optimal_level_0_new_qmc}, respectively, potentially leading to a finer model as level 0. This further distinguishes the adaptive scheme from the uniform approach.
\end{remark}

\begin{figure}[ht]
	\begin{subfigure}{0.48\textwidth}
		\includegraphics[width=\textwidth]{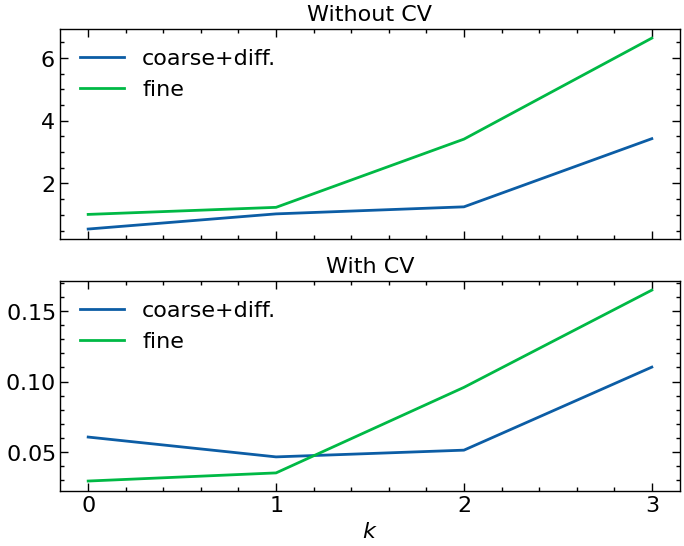}
		\caption{Uniform, MC. Top~\eqref{eq:optimal_level_0_old}, bottom~\eqref{eq:optimal_level_0_new_mc}.}	
	\end{subfigure}
	\begin{subfigure}{0.48\textwidth}
		\includegraphics[width=\textwidth]{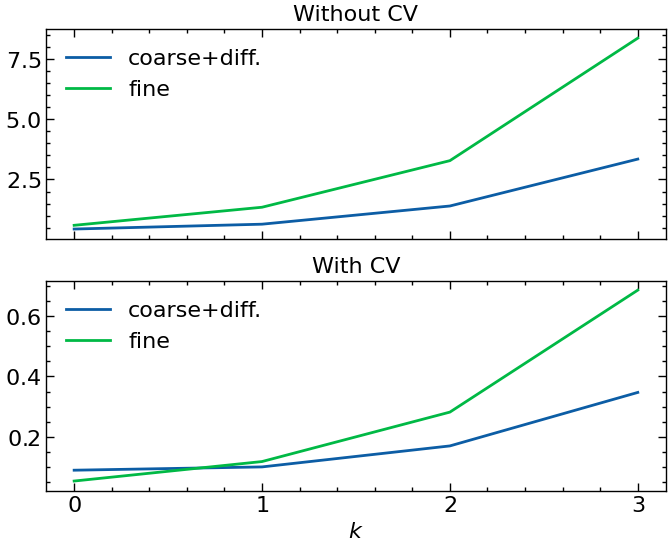}
		\caption{Uniform, QMC. Top~\eqref{eq:optimal_level_0_old_qmc}, bottom~\eqref{eq:optimal_level_0_new_qmc}.}	
	\end{subfigure}
	\\
	\begin{subfigure}{0.48\textwidth}
		\includegraphics[width=\textwidth]{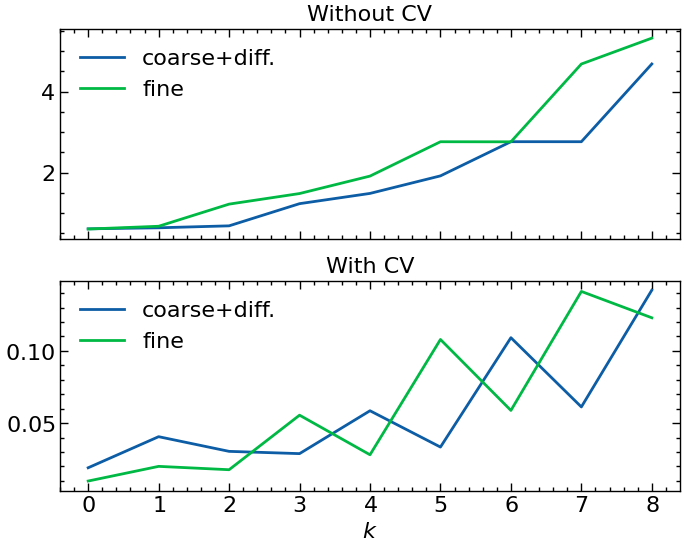}
		\caption{Adaptive, MC. Top~\eqref{eq:optimal_level_0_old}, bottom~\eqref{eq:optimal_level_0_new_mc}.}	
	\end{subfigure}
	\begin{subfigure}{0.48\textwidth}
		\includegraphics[width=\textwidth]{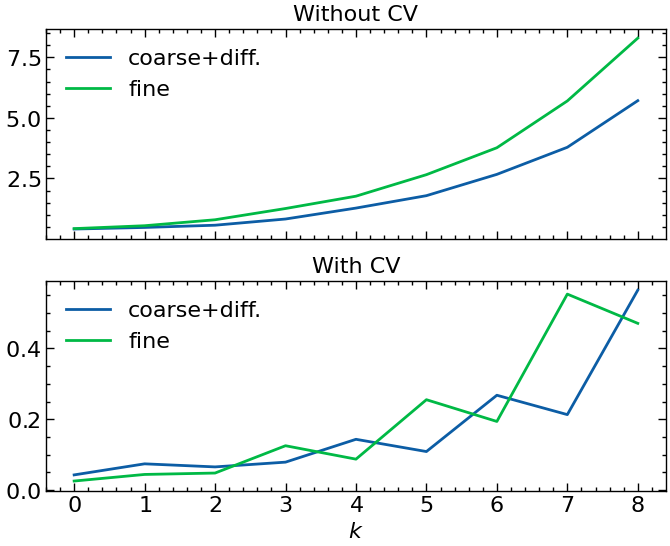}
		\caption{Adaptive, QMC. Top~\eqref{eq:optimal_level_0_old_qmc}, bottom~\eqref{eq:optimal_level_0_new_qmc}.}	
	\end{subfigure}
	\caption{Comparison of cost contributions between the coarse plus difference terms and the fine terms for Uniform and Adaptive meshes with MC and QMC. The ``coarse+diff.'' and ``fine'' corresponds to the left hand side and right hand side of the various inequalities referred in the subcaptions. {The optimal initial mesh is determined by the smallest index $k$ where the coarse plus difference term is smaller than the fine term.}}
	\label{fig:coarse_plus_difference_vs_fine}
\end{figure}

Figure~\ref{fig:coarse_plus_difference_vs_fine} compares the cost contributions between the coarse plus difference terms (the left hand side of inequalities~\eqref{eq:optimal_level_0_old}-\eqref{eq:optimal_level_0_new_qmc}) and the fine terms (the right hand side of inequalities~\eqref{eq:optimal_level_0_old}-\eqref{eq:optimal_level_0_new_qmc}) for the numerical example in the next section across uniform and adaptive meshes using both MC and QMC methods. The adaptive meshes follow the design in~\cite{adMLMC_our}. Across all scenarios (uniform/adaptive meshes with MC/QMC), the CV leads to a finer optimal mesh at the level 0.

\revRtwo{In our study we explore two families of level-0 control variates. The first is a \emph{truncation-based CV}: it is built from a truncated Karhunen--Loève (KL) expansion of the random field and admits two practical variants (``importance-ranked truncation'' (type-I) vs.\ ``natural-order truncation'' (type-II)). The second is an \emph{SVD-based CV}: it is derived from the singular value decomposition of an operator matrix and in general achieves stronger variance reduction in our experiments. }


For the \emph{truncation-based CV}, the control variate on level 0, denoted as $\tilde{Q}_{0}^t$, is the quantity of interest (QoI) obtained from the coefficient $a_t$, which includes only $t \ll s$ terms of the expansion:
\begin{equation}
a_t = \exp\left( \sum_{j=1}^t y_{J(j)} \psi_{J(j)} \right),
\end{equation}
where $J$ is a permutation of the indices $1, 2, \dotsc, s$, $y_j$ and $\psi_j$ are defined in~\eqref{eq:a_def}. This approach requires selecting $t$ important indices from the total $s$ indices. Methods such as the Sobol' indices can be employed to analyze the sensitivity of QoIs with respect to these indices. However, this method involves calculating the conditional expectation for each index, which can significantly increase the computational cost.

\begin{figure}[ht]
	\centering
	\includegraphics[width = 0.495\textwidth]{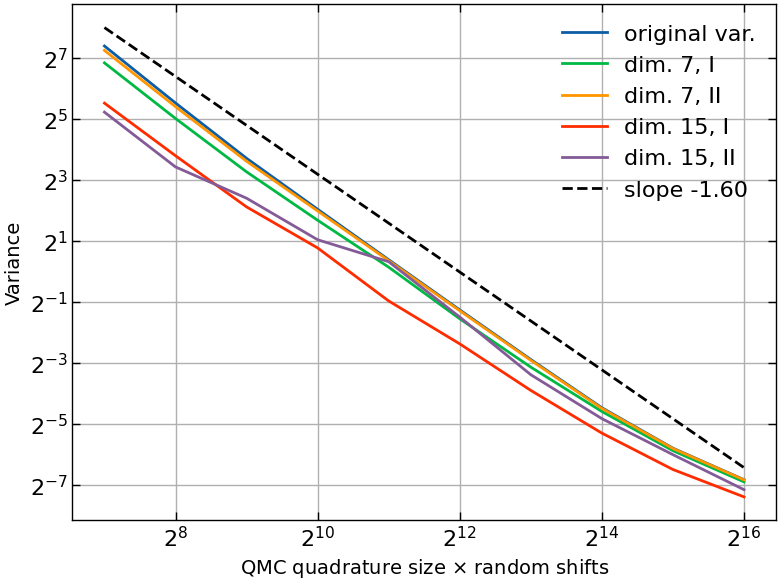}
	\includegraphics[width = 0.495\textwidth]{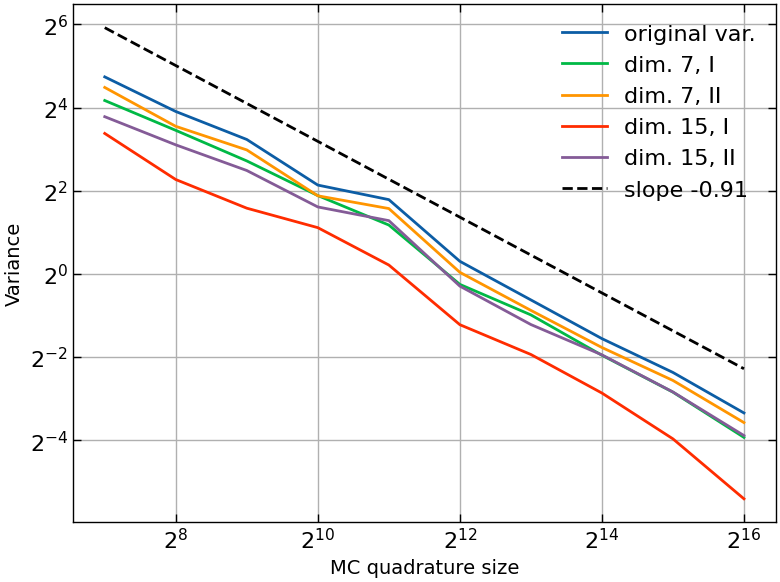}
	\caption{Variance of RQMC estimator (left) and MC estimator (right) for $\var{\tilde{Q}_0}$ (original var.) and the control variate $\var{\tilde{Q}_0 - \tilde{Q}_0^t}$ with dimension $t = 7, 15$ for both type-I, type-II CVs. }
	\label{fig:single_level_var_k_mc}
\end{figure}

To mitigate the costs associated with generating new data, we leverage the importance sampling technique discussed in Section~\ref{sec:theory_is}. Employing this method, we select the first $t$ dimensions characterized by the largest values of $\bar{\bm{\alpha}}_j$, as defined in~\eqref{eq:bar_alpha_is_optimization}. This selection criterion is informed by the model's behavior, where the upper bound of the QoI blows up at the boundaries. Larger values of $\bar{\bm{\alpha}}_j$ lead to a greater reduction to the boundary singularities, making it advantageous to apply importance sampling. We refer to this kind of low-dimensional CV as the \revRtwo{\emph{importance-ranked truncation CV}} (type-I) ; see Figure~\ref{fig:single_level_var_k_mc}. Selecting the first $t$ dimensions in the natural order of the series expansion in~\eqref{eq:a_def} is referred to as the \revRtwo{\emph{natural-order truncation CV}} (type-II). The variance of the RQMC estimators $\var{I_N(\tilde{Q}_0);\bm{\Delta}}$, $\var{I_N(\tilde{Q}_0 - \tilde{Q}_0^t);\bm{\Delta}}$ and MC estimators $\var{I_N(\tilde{Q}_0)}$, $\var{I_N(\tilde{Q}_0 - \tilde{Q}_0^t)}$ for the numerical example we consider in the next section is shown in Figure~\ref{fig:single_level_var_k_mc}.

In Figure~\ref{fig:single_level_var_k_mc}, we notice that {both truncation-based CVs reduce} the variance of the MC and RQMC estimators. The variance reduction is larger with more dimensions included in the CV. {The importance-ranked truncation CV, with the dimensions selected based on the criterion $\bar{\bm{\alpha}}_j$, is more effective than the natural-order truncation CV.}

If we hypothesize that sampling $\tilde{Q}_0 - \tilde{Q}_0^t$ incurs twice the cost of sampling $\tilde{Q}_0$, then for the control variate approach to be efficient, we require:
\begin{equation*}
	\var{I_N (\tilde{Q}_0 - \tilde{Q}_0^t)  } \leq \frac{1}{2} \var{I_N (\tilde{Q}_0)}	
\end{equation*}
and
\begin{equation*}
	\var{I_N (\tilde{Q}_0 - \tilde{Q}_0^t;\bm{\Delta})  } \leq \frac{1}{2} \var{I_N (\tilde{Q}_0;\bm{\Delta})}.
\end{equation*}
However as depicted in Figure~\ref{fig:single_level_var_k_mc}, the variance reduction is not \revRone{as large as} desired.

One more observation from Figure~\ref{fig:single_level_var_k_mc} is the different non-asymptotic behavior of the MC and RQMC estimators. To illustrate this phenomenon, we first establish lower bounds for the MC and RQMC estimators: 
\begin{equation}
\sqrt{\var{  I_N(\tilde{Q}_0 - \tilde{Q}_0^t)  } } \geq \sqrt{\var{  I_N(\tilde{Q}_0)  } } - \sqrt{\var{   I_N(\tilde{Q}_0^t) } },
\end{equation}
and for the RQMC setting:
\begin{equation}
	\sqrt{\textrm{Var} { \left[ I_N (\tilde{Q}_0 - \tilde{Q}_0^t; \bm{\Delta}) \right] } } \geq \sqrt{ \textrm{Var}  \left[ { {I}_{N} (\tilde{Q}_0;\bm{\Delta} ) } \right] } - \sqrt{\textrm{Var} \left[ { {I}_{N} ( \tilde{Q}_0^t;\bm{\Delta}) } \right] }. 
\end{equation}
As noted in~\cite{liu2023nonasymptotic}, the variance of the RQMC estimator for the truncated model, ${\textrm{Var} \left[ { {I}_{N} ( \tilde{Q}_0^t;\bm{\Delta}) } \right] }$ converges faster in the non-asymptotic regime than that of the original model ${\textrm{Var} \left[ { {I}_{N} ( \tilde{Q}_0;\bm{\Delta}) } \right] }$. This enhanced rate is attributed to the reduced nominal dimension, $t < s$, resulting in a reduced effective dimension. As the number of QMC quadrature points $N$ increases, the variance difference $\sqrt{\textrm{Var} { \left[ I_N (\tilde{Q}_0 - \tilde{Q}_0^t; \bm{\Delta}) \right] } }$ approaches $\sqrt{ \textrm{Var} \left[ { {I}_{N} (\tilde{Q}_0;\bm{\Delta} ) } \right] }$. 

The above considerations prompt us to explore an alternative CV strategy for effective variance reduction. Specifically, on level 0, we analyze the series expansion described in~\eqref{eq:a_def}. In this formulation, the logarithm of the coefficient values at the FEM quadrature points, denoted as $\bm{a}_0 \in \mathbb{R}^{q_0}$ are linear functionals of the $s$-dimensional normal random variable $\mathbf{y}_s$. Here, $q_0$ represents the number of FEM quadrature points on the mesh corresponding to level 0. The relationship is expressed as $\log(\bm{a}_0) = \bm{A}_0 \mathbf{y}_s$, where $\bm{A}_0 \in \mathbb{R}^{q_0 \times s}$. The QoI has a nonlinear dependence on these values at the quadratures, formulated as $\tilde{Q} = \tilde{Q}(\bm{A}_0 \mathbf{y}_s)$. Instead of truncating $\mathbf{y}_s$, as in the previous approach, we aim to develop a low-rank approximation of $\bm{A}_0$ for a low-dimensional representation of the CV for $\tilde{Q}$. 

We consider the singular value decomposition (SVD) of $\mat{A}_0$, given by
\begin{equation}
	\mat{A}_0 = \mat{U} \mat{\Sigma} \mat{V}^{T},
\end{equation}
with $\mat{U} \in \mathbb{R}^{q_0\times q_0}$, $\mat{\Sigma} \in \mathbb{R}^{q_0\times s}$ and $\mat{V}^{T} \in \mathbb{R}^{s\times s}$, where $q_0$ is the FEM quadrature size on mesh $0$. For $k \ll q,s$, a rank-$k$ approximation $\mat{A}_0^k$ is given by
\begin{equation}
	\mat{A}_0^k = \mat{U}_k \mat{\Sigma}_k \mat{V}^{T}_k,
\end{equation}
where $\mat{U}_k \in \mathbb{R}^{q_0 \times k}$, $\mat{\Sigma}_k \in \mathbb{R}^{k\times k}$ and $\mat{V}^{T}_k \in \mathbb{R}^{k\times s}$. The CV $\tilde{Q}^{\textrm{CV}}$ is then defined as: 
\begin{equation}
	\begin{split}
		\tilde{Q}^{\textrm{CV}} &= \tilde{Q}(\mat{A}_0^k \mathbf{y}_s )\\
		&= \tilde{Q}(\mat{U}_k \mat{\Sigma}_k \mat{V}^{T}_k \mathbf{y}_s ).
	\end{split}
\end{equation}
Since $\mathbf{y}_s \in \mathcal{N}(0, \mat{I}_s)$, we have that $\mat{V}^{T}_k \mathbf{y_s} \in \mathcal{N}(0, \mat{I}_k)$, due to the orthogonality of $\mat{V}^{T}_k$ ($\mat{V}^{T}_k \mat{V}_k = \mat{I}_k $). This lower-dimensional representation allows for efficient computation of $\mathbb{E}[\tilde{Q}^{\textrm{CV}}]$ using high-order quadrature methods tailored to the $k$-dimensional Gaussian measure, such as stochastic collocation. This approach significantly reduces the complexity of the calculations, especially given that $k \ll s$.

\begin{figure}[ht]
	\centering
	\begin{minipage}{0.52\linewidth}
		\begin{subfigure}{\linewidth}
		\includegraphics[width = \textwidth]{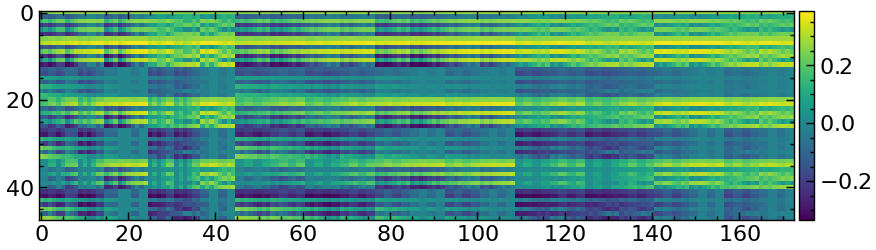}
		\caption{The matrix $\mat{A}_0$}\label{subfig:key-a}
		\end{subfigure}\\
		\begin{subfigure}{\linewidth}
		\includegraphics[width = \textwidth]{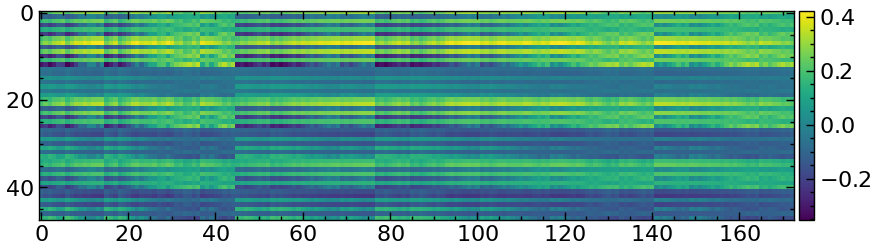}
		\caption{Rank-2 approximation}\label{subfig:key-b}
		\end{subfigure}\\
	\end{minipage}
	\begin{minipage}{0.4\linewidth}
		\begin{subfigure}{\linewidth}
			\includegraphics[width=\linewidth]{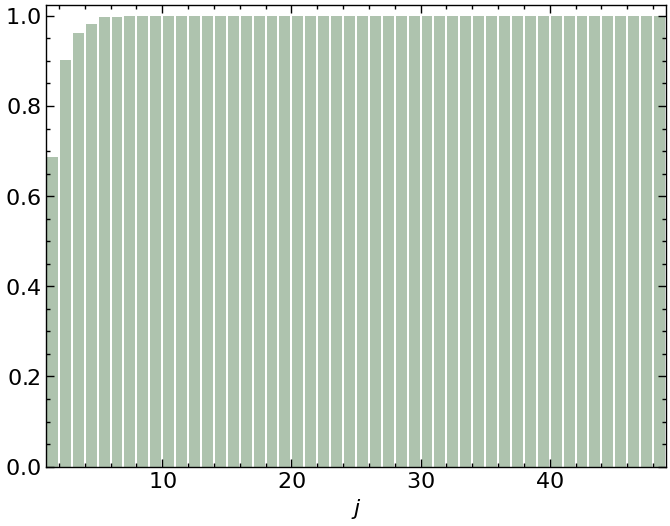}
			\caption{The ratio ${\sum_{i=1}^j \vartheta^2}/{\sum_{i=1}^{s} \vartheta^2}$}\label{subfig:key-d}
		  \end{subfigure}
	\end{minipage}
	\hfill
	\caption{{The matrix $\mat{A}_0$, its rank-2 approximation and the cumulative energy ratio (${\sum_{i=1}^j \vartheta^2}/{\sum_{i=1}^{s} \vartheta^2}$). The cumulative energy ratio at $j=2$ is around 0.9. }}
	\label{fig:matrix_a0}
\end{figure}

Figure~\ref{fig:matrix_a0} displays the matrix $\mat{A}_0$, its rank-2 SVD approximations (the color scale indicating the magnitude of matrix values) and the cumulative energy ratio of the singular values (${\sum_{i=1}^j \vartheta_i^2}/{\sum_{i=1}^{s} \vartheta_i^2}$). As shown in Figure~\ref{subfig:key-d}, the $L^2$ energy of the matrix $\mat{A}_0$ is predominantly captured by the initial singular values. Notice that the first two singular value accounts for approximately 90\% of the energy, and the rank-2 approximation is almost indistinguishable from the original matrix $\mat{A}_0$. We therefore use rank-2 SVD approximation for the CV $\tilde{Q}^{\textrm{CV}}$ in the numerical example.

\begin{figure}[ht]
	\includegraphics[width = 0.48\textwidth]{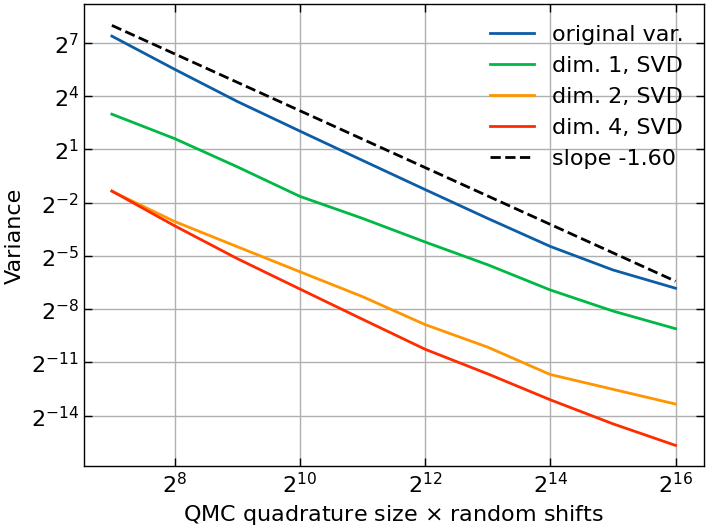}
	\includegraphics[width = 0.48\textwidth]{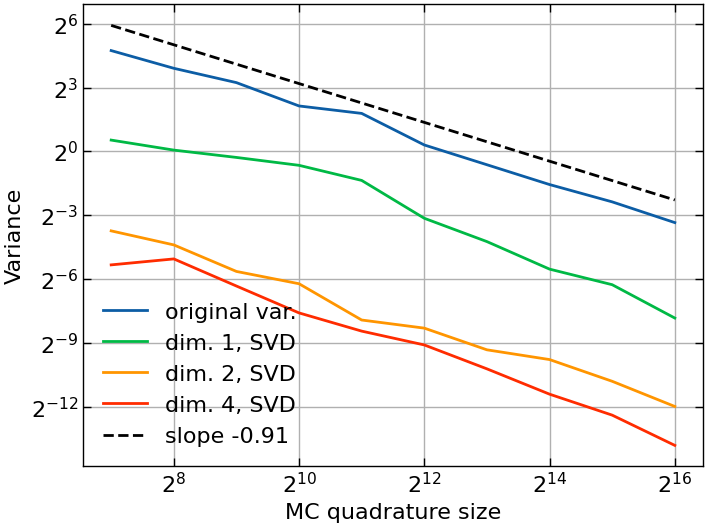}
	\hfill
	\caption{Variance of the SVD-based control variate RQMC estimator (left) $\sqrt{\textrm{Var} { \left[ I_N (\tilde{Q}_0 - \tilde{Q}_0^t; \bm{\Delta}) \right] } }$ and MC estimator (right) $\sqrt{\textrm{Var} { \left[ I_N (\tilde{Q}_0 - \tilde{Q}_0^t) \right] } }$.}
	\label{fig:svd_cv}
\end{figure}

Figure~\ref{fig:svd_cv} presents the variance of the SVD-based control variate for both the RQMC estimator and the MC estimator. We observe a significant variance reduction through the use of the SVD-based control variate. 

\section{Numerical Results}
\label{sec:numex}

\revRtwo{\subsection{Numerical workflow}\label{sec:numex_workflow}
Before presenting the results, we record the workflow and the discretization choices that apply throughout this section, so that the complexity rates reported below can be interpreted against a concrete setup.
\begin{description}
\item[Geometry and QoI.] 2-D slit domain $\mathcal{D}=[-1,1]\times[-1,0]$ with Dirichlet boundary $\partial\mathcal{D}_1=\partial\mathcal{D}\setminus([-1,0]\times\{0\})$ and QoI given by~\eqref{eq:QoI_general} (a convolution of the indicator of the target subregion $[0.25,0.5]\times[-0.5,-0.25]$ with a fixed Gaussian kernel).
\item[Coefficient.] Lognormal~\eqref{eq:a_def} with a truncated Mat\'ern series expansion (smoothness $\nu=4.5$, correlation length $\varrho=1.0$), truncated to $s=49$ Gaussian random variables ranked by descending eigenvalue of the covariance operator.
\item[FEM discretization.] Piecewise-bilinear Lagrange ($Q_1$) elements on structured quadrilateral meshes with hanging nodes, i.e.\ polynomial degree $p=1$ in each coordinate direction. The adaptive-mesh hierarchy is produced by the goal-oriented $h$-adaptive solver of~\cite{adFEM_our}, with stopping tolerance $\tol_\ell=2^{-\ell-2}$ per level unless stated otherwise.
\item[Quadrature.] MC samples use standard pseudorandom points; QMC uses a scrambled Sobol' sequence. The number of randomizations $R$ is \emph{not} uniform across the section: we use $R=64$ for the level-wise variance-decay plots (Figures~\ref{fig:ex_2_qmc_convergence_is} and \ref{fig:mlqmc_u_a_cv_cost}) because reducing the MC error on the variance estimator is the bottleneck for those plots, and we use $R=8$ for the RMSE verification in Figure~\ref{fig:mlqmc_cv_error_verification} because that plot already averages over 100 independent outer realizations and adding randomizations would inflate wall-clock time without changing the qualitative conclusion.
\item[Variance reduction.] Importance sampling uses a Gaussian proposal $\varphi_{\bm\alpha}$ constructed from the optimization of Section~\ref{sec:theory_is}, with an offline training sample of size $n$ drawn from the target standard Gaussian density $\varphi$ before the production runs. The specific value of~$n$ used in each experiment is fixed in advance and is kept small so that its offline cost is negligible relative to the MLQMC work reported in Figures~\ref{fig:mlqmc_u_a_cv_cost} and~\ref{fig:mlqmc_cv_error_verification}. The level-0 control variates use the truncation-based and SVD-based variants.
\item[Reference solution.] The reference expectation used to compute RMSE in Figure~\ref{fig:mlqmc_cv_error_verification} is obtained by a high-accuracy MLQMC-CV run at $\tol=2^{-12}$.
\end{description}}

This section presents numerical results using a random coefficient modeled as a series expansion inspired by the Mat\'ern covariance model. \revRtwo{The numerical example uses} the same PDE model, domain, and quantity of interest (QoI) as previously described in~\cite{adMLMC_our}. For completeness, we state the PDE model and the QoI below. Recall the PDE model~\eqref{eq:pde_general_1}:
\begin{subequations}
	\begin{align}
	-\nabla \cdot \left( a(\mathbf x; \mathbf{y}) \nabla u(\mathbf x; \mathbf{y}) \right) 
	&= f(\mathbf x) &&\text{for $\mathbf x  \in \mathcal D$,} \nonumber
	\\
	u(\mathbf x; \mathbf{y}) &=  0 &&\text{for $\mathbf x \in \partial\mathcal{D}_1$},\nonumber
	\\
	\partial_n u(\mathbf x; \mathbf{y}) &=  0 &&\text{for $\mathbf x \in \partial\mathcal{D} -\partial\mathcal{D}_1$},\nonumber
	\end{align}
\end{subequations}
where $\mathcal{D} = [-1, 1]\times [-1, 0] \subset \mathbb{R}^2$, $\partial \mathcal{D}_1 = \partial \mathcal{D} \setminus \left( [-1, 0] \times \{0\} \right)$, and $f(\mathbf x) = 1$. The QoI is given by the following:
\begin{equation}
	\label{eq:QoI_general}
	Q = \int_{\mathcal{D}} u(\mathbf x; \mathbf{y})  (\mathbbm{1}_{[0.25, 0.5]\times [-0.5, -0.25]} * \varphi ) (\mathbf{x}) \, \mathrm{d} \mathbf x.
\end{equation}
where $*$ denotes convolution and $\varphi$ is a 2-d Gaussian kernel, given by
\begin{equation*}
	\varphi(\mathbf{x}) = \frac{8}{\pi} \exp\left(-8 \mathbf{x}^T  \mathbf{x} \right).
\end{equation*}

Recall the coefficient model~\eqref{eq:a_def}:
\begin{equation}
	a = \exp\left(\sum_{j=1}^{s} \psi_j \mathbf{y}_j \right)
\end{equation}
where $\psi_j = \sqrt{\lambda_j} \theta_j$, with $\lambda_j, \theta_j$ denoting the eigenvalues and orthonormal eigenfunctions of the Mat\'ern covariance kernel with smoothness parameter $\nu = 4.5$ and correlation length $\varrho= 1.0$. The terms $\psi_j$ are ranked in descending order of $\lambda_j$, and we select $s = 49$ in the expansion~\eqref{eq:a_def}.

\begin{figure}[ht]
	\begin{subfigure}{0.30\linewidth}
	\includegraphics[width = \textwidth]{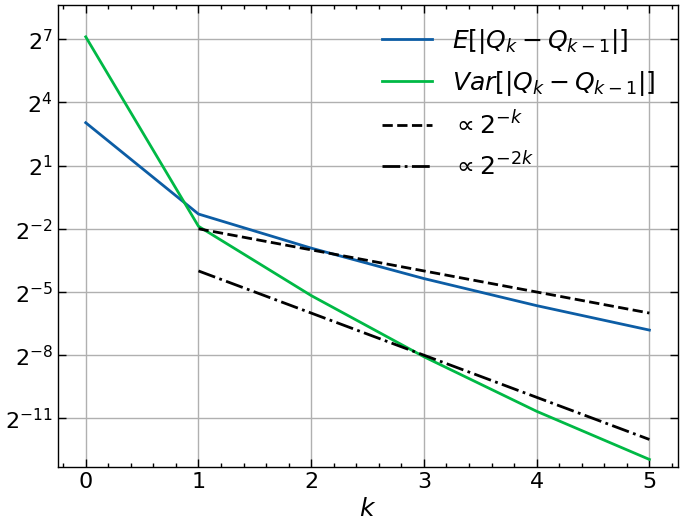}
	\caption{Uniform Mesh}
	\end{subfigure}
	\begin{subfigure}{0.30\linewidth}
	\includegraphics[width = \textwidth]{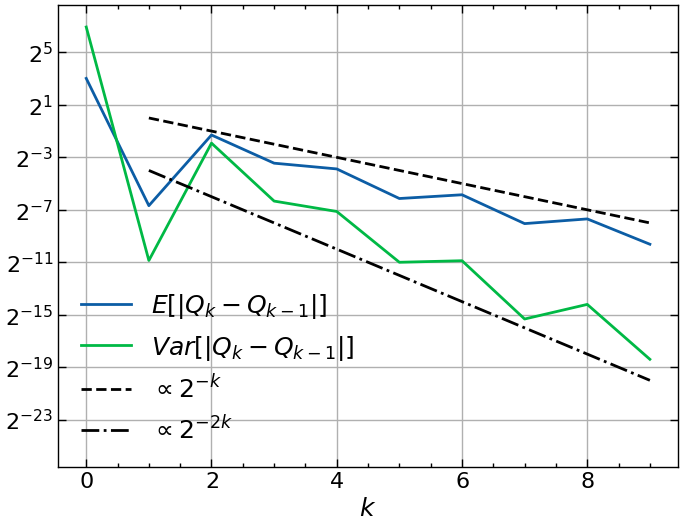}
	\caption{Adaptive Mesh}
	\end{subfigure}
	\begin{subfigure}{0.30\linewidth}
	\includegraphics[width = \textwidth]{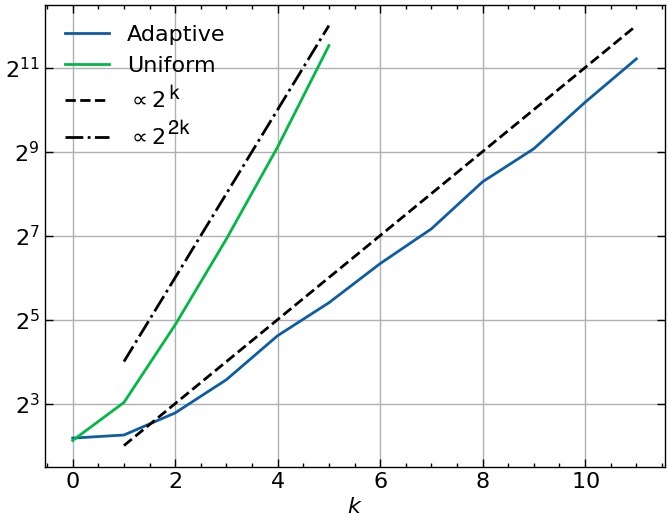}
	\caption{Cost}
	\end{subfigure}
	\hfill
	\caption{Decay of mean and variance against level $k$ for uniform Meshes (left) and adaptive Meshes (middle). Growth of wall-clock time against $k$ (right).}
	\label{fig:single_level_k_mc}
\end{figure}

\begin{figure}[ht]
	\centering
	\hfill
	\includegraphics[width = 0.495\textwidth]{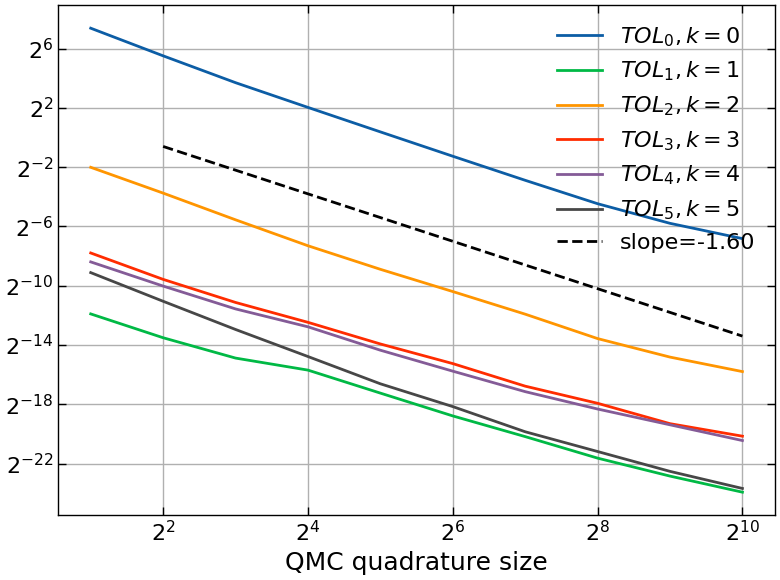}
	\includegraphics[width = 0.495\textwidth]{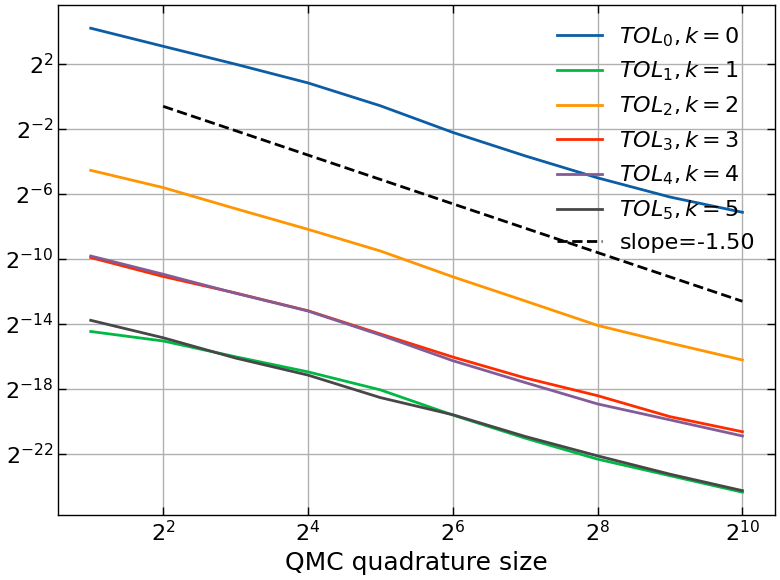}
	\caption{Estimated variance $\textrm{Var}  \left[{I_{N_k} (\Delta \tilde{Q}_{k}; \bm{\Delta}) }\right] $ for $k = 0, 1, 2, \dotsc, 5$ with $\tol_\ell = 2^{-\ell - 2}$, $\ell = 0, 1, 2, \dotsc, 5$, with 64 random shifts. Results are shown without IS (left) and with IS (right). \revRtwo{The importance sampling does not materially improve the convergence rate in this range, but produces a per-level variance reduction visible across all terms plotted.}}
	\label{fig:ex_2_qmc_convergence_is}
\end{figure}

To analyze the complexity of the multilevel estimators, we first examine the decay of the mean and variance against level $k$ for uniform and adaptive meshes. Figure~\ref{fig:single_level_k_mc} plots the mean, variance and computational cost for both mesh types, along with the fitted the convergence rates. Adaptive meshes achieve the same mean and variance convergence rates as uniform meshes but with reduced computational complexity. 

Figure~\ref{fig:ex_2_qmc_convergence_is} plots the convergence of the estimated variance $\textrm{Var}  \left[{I_{N} ( \Delta \tilde{Q}_{k}; \bm{\Delta}) }\right]$ for different levels $k$ against the QMC quadrature size $N$. Results are shown without IS (left) and with IS (right), where the type of IS is detailed in Section~\ref{sec:theory_is}. The variance is estimated with randomizations:
\begin{equation*}
	\textrm{Var}  \left[{I_{N} (\Delta \tilde{Q}_{k}; \bm{\Delta}) }\right] \approx \frac{1}{R - 1} \sum_{r=1}^R \left( {I_{N} (  \Delta \tilde{Q}_{k}; \bm{\Delta}_r) } - \frac{1}{R} \sum_{r=1}^R I_{N} ( \Delta \tilde{Q}_{k}; \bm{\Delta}_r) \right)^2. 
\end{equation*}
As analyzed in~\cite{liu2023nonasymptotic}, this type of RPDE model satisfies a certain boundary growth condition, leading to an RQMC variance convergence rate of $\mathcal{O}(n^{-2+\epsilon})$ for any $\epsilon > 0$. In the plotted pre-asymptotic regime, the observed convergence rate is found to be $\lambda = 1.6$. Although the IS does not improve the convergence rate in the observed range, it reduces the variance by a factor of approximately 2 across all terms plotted; As noted in~\cite{liu2023nonasymptotic}, this type of IS becomes more effective with the increased coefficient variability. 

\begin{figure}[ht]
	\centering
	\hfill
	\includegraphics[width = 0.495\textwidth]{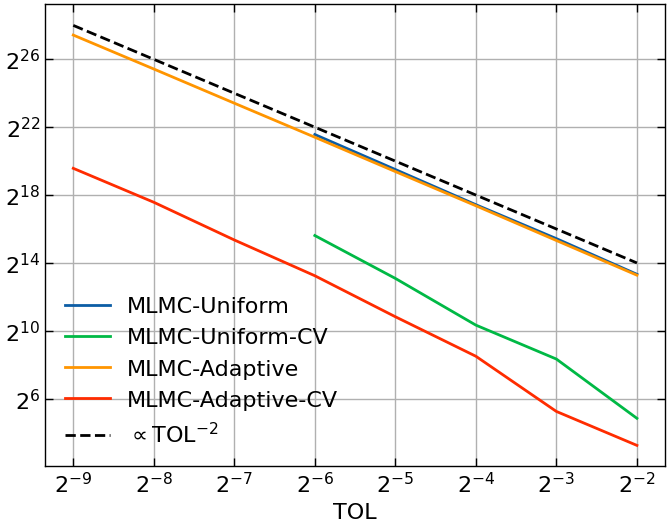}
	\includegraphics[width = 0.495\textwidth]{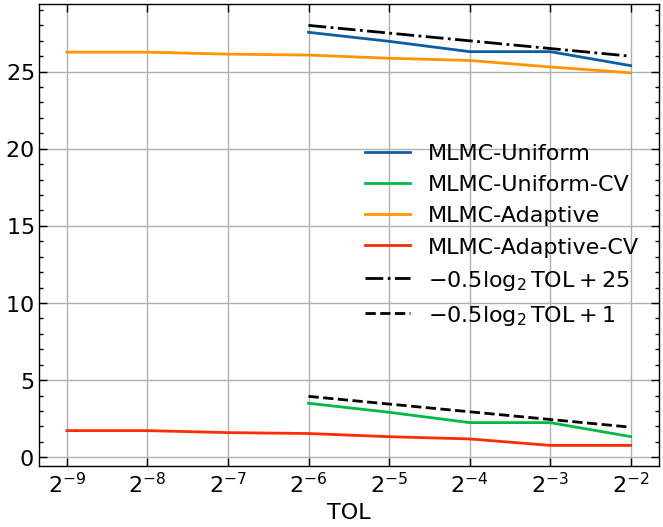}\\
	\caption{The estimated MLMC cost (left) and the factor $\left( \sum_{k=0}^{L(\tol)} 
	\sqrt{V_k C_k} \right)^{2}$ (right). The advantage of adaptive mesh is clear when the CV is applied. }
	\label{fig:mlmc_u_a_cv_cost}
\end{figure}
\begin{figure}[ht]
	\centering
	\hfill
	\includegraphics[width = 0.495\textwidth]{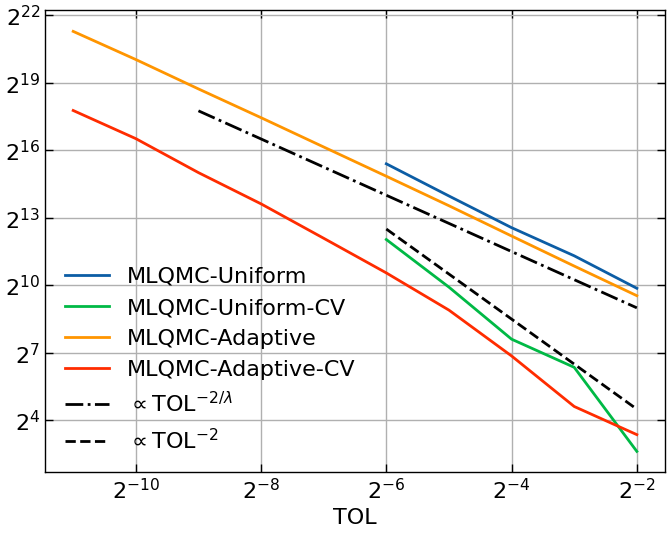}
	\includegraphics[width = 0.495\textwidth]{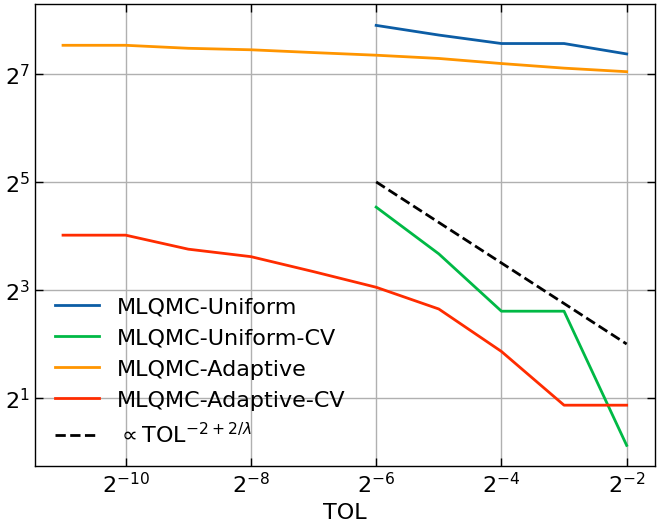}\\
	\caption{The estimated MLQMC cost (left) and the factor $ \left( \sum_{k=0}^{L(\tol)} C_k
	\left(\frac{V_k^{\mathrm{Q}}}{C_k}\right)^{\frac{1}{\lambda + 1}} \right)^{\frac{\lambda + 1}{\lambda}} R^{1- \frac{1}{\lambda}} $ (right) with $\lambda = 1.6$. \revRtwo{In the parameter range explored, the adaptive meshes improve the asymptotic work scaling of the MLQMC method with and without the CV.} }
	\label{fig:mlqmc_u_a_cv_cost}
\end{figure}
\begin{figure}[htbp]
	\centering
	\hfill
	\includegraphics[width = 0.495\textwidth]{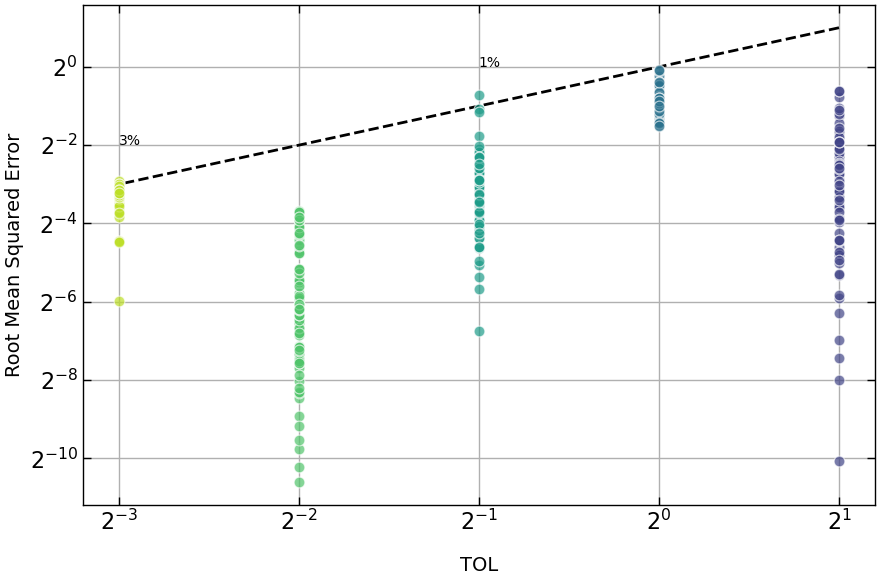}
	\includegraphics[width = 0.495\textwidth]{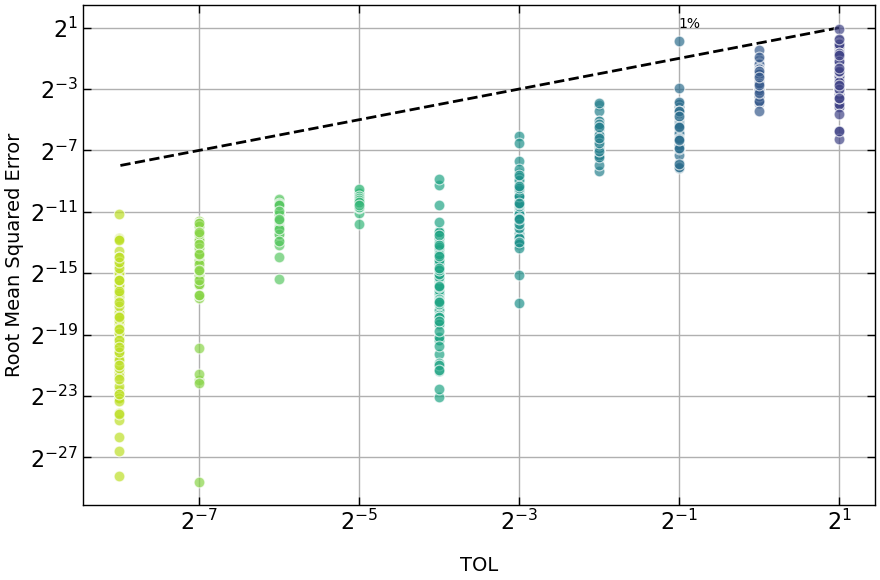}
	\caption{\revRtwo{Error verification for MLQMC (left) and MLQMC with CV (right) on adaptive meshes. The reference solution is obtained by a MLQMC-CV estimator at $\tol = 2^{-12}$ (same reference used in both panels). The dashed lines indicate where the root mean squared error equals $\tol$. The percentage of outliers is marked in the figure. }}
	\label{fig:mlqmc_cv_error_verification}
\end{figure}

Recall that the computational work of MLMC is given by
\begin{equation*}
	\tol^{-2} \left( \sum_{k = 0}^{L(\tol)} \sqrt{V_k C_k} \right)^2.
\end{equation*}
Figure~\ref{fig:mlmc_u_a_cv_cost} plots both the estimated MLMC cost and the factor $\left( \sum_{k=0}^{L(\tol)} \sqrt{V_k C_k} \right)^{2}$. Without the CV, MLMC costs for uniform and adaptive meshes are \revAE{almost} indistinguishable due to level 0 dominating the computational cost in both cases. Implementing the CV alongside the optimal initial mesh reveals the advantages of adaptive meshes: CV on the initial level significantly reduces the cost and reveals the improved convergence rate compared to the uniform meshes.

The right plot shows that the factor $\left( \sum_{k=0}^{L(\tol)} \sqrt{V_k  C_k} \right)^{2}$ grows linearly with $\log (\mathrm{TOL})$ for uniform meshes, both with and without the CV. However, for adaptive meshes, the term $\left( \sum_{k=0}^{L(\tol)} \sqrt{V_k  C_k} \right)^{2}$ grows much more slowly and appears to stabilize around a constant as $\tol$ decreases. Table~\ref{tab:mlmc_cost_vlcl} provides the numerical estimates of $\sqrt{V_k C_k}$. \revAE{Note that} when using the CV and the optimal initial mesh, the first uniform mesh and the first three adaptive meshes are discarded from the mesh hierarchy.

\begin{table}[ht]
	\centering
	\begin{tabular}{lllllll}
	\hline
	 $k$ & 0 & 1 & 2 & 3 & 4 & 5\\
	 \hline
	 MLMC-U. & 24.54 & 1.48 & 0.90 & 0.67 & 0.58 & 0.61 \\
	 MLMC-U.-CV & - & 1.34 & 0.90 & 0.67 & 0.58 & 0.61 \\
	 \hline
	\end{tabular}
	\begin{tabular}{lllllllllll}
		\hline
		 $k$ & 0 & 1 & 2 & 3 & 4 & 5 & 6 & 7 & 8 & 9\\
		 \hline
		 MLMC-A. & 23.52 & 0.05 & 1.35 & 0.38 & 0.42 & 0.14 & 0.21 & 0.06 & 0.13 & 0.04\\
		 MLMC-A.-CV & - & - & - & 0.77 & 0.42 & 0.14 & 0.21 & 0.06 & 0.13 & 0.04 \\
		 \hline
	\end{tabular}
	\caption{Cost contribution from each level $\sqrt{V_k C_k}$. }
	\label{tab:mlmc_cost_vlcl}
\end{table}

Next, we present the results for MLQMC. Figure~\ref{fig:mlqmc_u_a_cv_cost} plots both the estimated MLQMC cost and the factor $\left( \sum_{k=0}^{L(\tol)} C_k
\left(\frac{V_k^{\mathrm{Q}}}{C_k}\right)^{\frac{1}{\lambda + 1}} \right)^{\frac{\lambda + 1}{\lambda}} R^{1- \frac{1}{\lambda}}$ for uniform and adaptive meshes, with and without CV. Similar to the MLMC case, computational cost for uniform and adaptive meshes are indistinguishable without CV, while the CV reveals the advantages of adaptive meshes. For adaptive meshes, the factor $ \left( \sum_{k=0}^{L(\tol)} C_k \left(\frac{V_k^{\mathrm{Q}}}{C_k}\right)^{\frac{1}{\lambda + 1}} \right)^{\frac{\lambda + 1}{\lambda}} R^{1- \frac{1}{\lambda}} $ shows no dependence on $\log (\textrm{TOL})$ for adaptive meshes. Table~\ref{tab:mlqmc_cost_vlcl} provides numerical estimates of the factor ${V_k^{\mathrm{Q}}}^{\frac{1}{1+\lambda}} {C_k}^{\frac{\lambda}{1+\lambda}}$ for each level $k$. \revRtwo{Just as for MLMC, for MLQMC when using the CV and the optimal initial mesh, the first uniform mesh and the first three adaptive meshes are discarded from the mesh hierarchy.}

\begin{table}[ht]
	\centering
	\begin{tabular}{lllllll}
	\hline
	 $k$ & 0 & 1 & 2 & 3 & 4 & 5\\
	 \hline
	 MLQMC-U. & 20.52 & 2.64 & 1.99 & 1.74 & 2.13 & 5.55 \\
	 MLQMC-U.-CV & - & 1.05 & 1.99 & 1.74 & 2.13 & 5.55 \\
	 \hline
	\end{tabular}
	\begin{tabular}{lllllllllll}
		\hline
		 $k$ & 0 & 1 & 2 & 3 & 4 & 5 & 6 & 7 & 8 & 9\\
		 \hline
		 MLQMC-A. & 18.11 & 0.11 & 1.92 & 0.57 & 0.76 & 0.88 & 0.58 & 0.48 & 0.52 & 0.28\\
		 MLQMC-A.-CV & - & - & - & 1.45 & 0.76 & 0.88 & 0.58 & 0.48 & 0.52 & 0.28\\
		 \hline
	\end{tabular}
	\caption{Cost contribution from each level ${V_k^{\mathrm{Q}}}^{\frac{1}{1+\lambda}} {C_k}^{\frac{\lambda}{1+\lambda}}$. }
	\label{tab:mlqmc_cost_vlcl}
\end{table}
\revRtwo{Finally we present the numerical verification results for the MLQMC and MLQMC with CV on adaptive meshes in Figure~\ref{fig:mlqmc_cv_error_verification}. Specifically, we plot the root mean squared error (RMSE) for 100 independent realizations of the MLQMC estimator, each based on $R = 8$ randomly-scrambled Sobol' sequences, against the prescribed tolerance $\tol$ (the $R=8$ choice is discussed in Section~\ref{sec:numex_workflow}). The maximum level $K$ is selected based on the bias error constraint~\eqref{eq:mlmc_bias_error}, while the number of samples $N_k$ on each level $k = 0, \dotsc, K$ is determined by the statistical error constraint~\eqref{eq:mlmc_statistical_error}. The reference expectation used to define the RMSE is obtained from a single high-accuracy MLQMC-CV run at $\tol = 2^{-12}$, and the same reference is used for both the left (MLQMC) and right (MLQMC-CV) plots so that the two panels are directly comparable. The dashed lines indicate where the RMSE equals $\tol$. Both MLQMC and MLQMC with CV achieve the desired accuracy with few outliers. With the CV, we are able to reach smaller $\tol$ values due to the reduced computational cost.}

This section demonstrates the effectiveness of the CV in reducing computational cost for both MLMC and MLQMC methods. While the CV preserves asymptotic complexity, it significantly reduces constant factors in the cost. Adaptive meshes become more efficient than uniform meshes when combined with the CV, and the \revAE{proposed adaptive} MLQMC method achieves lower costs than \revAE{the corresponding adaptive} MLMC due to QMC's faster variance decay.
\section{Conclusions}
\label{sec:concl}
In this work, we propose an adaptive MLQMC method for the goal-oriented approximation of a linear elliptic PDE with fixed geometric singularities and a lognormal random diffusivity coefficient. \revRtwo{The target problem class combines four analytical difficulties identified in the introduction: (D1) a geometric boundary singularity, (D2) a lognormal coefficient field without a deterministic positive lower bound, (D3) sample-dependent adaptive-mesh stopping that introduces parameter-space discontinuities, and (D4) infinitely many discontinuity locations that preclude classical pre-integration smoothing.}

In~\cite{adMLMC_our}, we introduced sample-dependent adaptivity within the MLMC algorithm to reduce computational costs in approximating the QoI. \revRtwo{That setting was well suited to MC sampling, since the stochastic mesh selection could be used directly without imposing additional regularity requirements on the parametric dependence of the integrand. For QMC methods, however, the same feature becomes a central difficulty: the sample-dependent stopping rule introduces discontinuities in parameter space, and in the present lognormal setting these occur at infinitely many a priori unknown locations. In this work, we do not attempt a direct QMC analogue of the earlier adaptive MLMC construction. Instead, we focus on deterministic adaptive hierarchies combined with level-wise randomized QMC, while retaining the geometric advantages of adaptivity.}

\revRtwo{The present work extends the earlier adaptive MLMC framework of~\cite{adMLMC_our} to the QMC setting by replacing the earlier sample-dependent construction with deterministic adaptive hierarchies and level-wise randomized QMC, and by complementing this with variance reduction to reduce the dominant coarse-level cost and thereby make effective use of the adaptive hierarchy on the benchmark studied here.}

\revRtwo{Our complexity discussion identifies the regimes in which adaptive hierarchies are expected to be advantageous by relating the level-wise bias decay and variance decay to the effective work-growth exponent. The numerical results on the 2-D slit benchmark show that adaptive and uniform hierarchies exhibit comparable decay in the mean and variance of the multilevel corrections, while the adaptive hierarchy exhibits reduced work growth. The randomized QMC experiments show a pre-asymptotic variance decay close to $N^{-1.6}$, and the resulting adaptive MLQMC cost, with bilinear Lagrange FEM ($Q_1$, polynomial degree $p=1$ on quadrilateral meshes), lies between $\mathcal{O}(\tol^{-1})$ and $\mathcal{O}(\tol^{-2})$, consistent with the analysis of~\eqref{eq:mlqmc_adaptive_2d_complexity} in the best regime of Table~\ref{tab:complexity_summary}.}

\revRtwo{We investigated two variance reduction approaches within the multilevel framework. The first approach employs importance sampling, leveraging the lognormal distribution of the diffusivity coefficient. In the parameter range studied, this gives a moderate variance reduction and primarily serves as a regularity-preserving change of measure. The second approach introduces a CV at level 0, which is particularly useful when level 0 dominates the total cost. We consider two constructions of the control variate: the first truncates the series expansion of the random field, which is less effective when the leading eigenvalues share similar magnitude; the second uses an SVD-based low-rank mapping from the input space to the logarithm of the random field, which is more effective in our numerical results. Moreover, incorporating the level-0 control variate can change the optimal choice of the initial mesh, thereby reducing the dominant coarse-level contribution and revealing the lower multilevel cost of the adaptive hierarchy.}

\revRtwo{Numerical results demonstrate that, on this 2-D slit benchmark and in the parameter range explored, adaptive MLQMC achieves the desired tolerance at a markedly lower computational cost than standard MLMC on the same mesh hierarchy; further reductions in computational cost are obtained when the level-0 control variate is included.}

\revRtwo{We emphasize three limitations of the present work. First, all numerical evidence is drawn from a single 2-D slit benchmark with a bilinear-Lagrange ($Q_1$, $p=1$) FEM discretization; the complexity rates in Table~\ref{tab:complexity_summary} are theoretical and would need to be re-verified on genuine higher dimensional problems or on coefficient fields with substantially different correlation structure. Second, the importance-sampling and level-0 control-variate constructions rely on an offline training phase whose cost is not charged against the reported work. Third, our model assumes a deterministic, geometric-driven reason for adaptivity. Problems that require stochastic remeshing are outside the scope of this paper. }

In future work, we plan several extensions to broaden and deepen our adaptive MLQMC framework. First, we will investigate fully adaptive control variates that evolve across levels; this can potentially change the weak and strong error convergence rates, as in Assumption~\ref{assumption:mc_complexity} and Assumption~\ref{assumption:qmc_rates}, \revRtwo{and thereby improve} the complexity. We also aim to tackle time-dependent PDEs; our adaptive MLQMC methods hold promise for applications in subsurface flow and transport modeling, uncertainty quantification in structural and petroleum engineering, and other relevant areas. At the same time, we recognize that very high or infinite stochastic dimensions pose challenges for geometry-driven adaptivity. \revRtwo{In such regimes, adaptivity is driven more by coefficient regularity,} geometric considerations become less influential, and alternative quadrature designs may prove more effective.
\section*{Acknowledgements}

This publication is based on work supported by the Alexander von Humboldt Foundation and the King Abdullah
University of Science and Technology (KAUST) office of sponsored research (OSR) under Award
No. OSR-2019-CRG8-4033. This work utilized the resources of the Supercomputing Laboratory at King Abdullah University of Science and Technology (KAUST) in Thuwal, Saudi Arabia.
We also acknowledge the use of the following open-source software
packages: \texttt{deal.II}~\cite{dealII92}.

\appendix
\section{Piecewise integration and measurability}
\label{sec:adaptivity_qmc_piecewise_integration}
This section explores the sample-adaptivity within the QMC framework. We introduce a piecewise integration and summation-by-parts formulation to compute the expectation of the QoI. \revRone{The present appendix is therefore exploratory: although conceptually natural, the numerical experiments reported below showed that the piecewise/summation-by-parts formulation did not deliver a significant improvement over the variance-reduction machinery of Section~\ref{sec:variance_reduction}, and it is \emph{not} part of the proposed method advocated in the main text. We include the derivation here for completeness, to document the approach, and to separate cleanly the two ideas that it contains: (i) the \emph{parameter-space discontinuity structure} of the sample-dependent adaptive FEM, which is the feature that motivates our variance-reduction design, and (ii) the \emph{infinite-sum representation} obtained by integrating analytically over $y_1$, which is what this appendix explores but which does not end up paying off numerically on the slit benchmark.}

For a given $\tol$ and $k \in \mathbb{N}_0$, define the region
\begin{align}
Y_{k, \tol} = \left\{ \mathbf{y} \in \mathbb{R}^s: \mathcal{K}(\mathbf{y}, \tol) = k \right\}.
\end{align}
where $\mathcal{K}(\mathbf{y}, \tol)$ is the mesh-selection function from~\eqref{eq:sample_dependent_stopping_discrete_deterministic}. \revRone{Writing
\begin{equation}
\label{eq:gj_definition}
  g_j(\mathbf y,\tol) \coloneqq \sum_K \rho_{j,K}(\mathbf y)\,h_{j,K}^{p+d} \; - \; \tol \, \frac{\int_{\mathcal D} \rho(\mathbf x;\mathbf y)^{d/(p+d)}\,d\mathbf x}{\int_{\mathcal D} \mathbb E[\rho^{d/(p+d)}]\,d\mathbf x}
\end{equation}
for the level-$j$ residual of the stopping criterion~\eqref{eq:sample_dependent_stopping_discrete_deterministic}, the level set admits the explicit representation
\begin{equation}
\label{eq:Yk_set_identity}
  Y_{k,\tol} = \Bigl\{ \mathbf y \in \mathbb R^s \; : \; g_0(\mathbf y,\tol) > 0, \; \dots, \; g_{k-1}(\mathbf y,\tol) > 0, \; g_k(\mathbf y,\tol) \le 0 \Bigr\}.
\end{equation}
Each $g_j(\cdot,\tol)\colon \mathbb R^s \to \mathbb R$ is continuous in~$\mathbf y$ (inherited from the continuity of the error-indicator terms $\rho_{j,K}$ and of the conditional normalizer), so $Y_{k,\tol}$ is a finite intersection of open and closed half-spaces of continuous functions and hence Borel-measurable. More generally, $\mathcal{K}(\cdot,\tol)\colon\mathbb{R}^s\to\mathbb{N}_0$ takes at most countably many values and is the first-hitting index of the sequence $\{g_j\}_{j\ge 0}$, so each $Y_{k,\tol}$ is Borel-measurable; this is all that is needed for the change-of-order-of-integration below to be rigorous. In particular, we do \emph{not} require the earlier assumption that each $Y_{k,\tol}$ be simply connected: the QoI $\bar{Q}_k$ is single-valued on each $Y_{k,\tol}$ regardless of its connected-component structure, so the identities~\eqref{eq:summation_Q}--\eqref{eq:summation_by_parts} below are valid as identities of Lebesgue integrals.}

Let $\bar{Q}(\mathbf{y}; \tol)$ denote the adaptive approximation of $Q$. Let $\bar{Q}_k (\mathbf{y})$ be the evaluation of $\bar{Q}(\mathbf{y})$ on mesh $k$. We write:
\begin{align*}
	\bar{Q}(\mathbf{y}; \tol) =  \bar{Q}_{\mathcal{K}(\mathbf{y},\tol)}(\mathbf{y}).
\end{align*}
In general the function $\bar{Q}(\mathbf{y}; \tol)$ is discontinuous in $\mathbf{y}$, \revRone{a consequence of mesh switching in the adaptive solver: each mesh corresponds to a distinct level set $Y_{k,\tol}$ of the stopping function, and the piecewise-smooth restrictions $\bar Q|_{Y_{k,\tol}}$ do not match across level-set boundaries. } Hence the expectation $\mathbb{E}[\bar{Q}(\mathbf{y}; \tol)]$ can be \revRone{decomposed into a sum of contributions from the individual level sets $Y_{k,\tol}$; since $\{Y_{k,\tol}\}_{k\ge 0}$ is a countable Borel partition of $\mathbb{R}^s$ and the integrand $\bar Q(\mathbf y;\tol)\,\varphi(\mathbf y)$ is non-negative, Fubini-Tonelli's theorem~\cite[Thm.~8.8]{Rudin1987} justifies the interchange of sum and integral below. Under the mild assumption that arbitrarily fine meshes are available in the solver, this decomposition can be} written as evaluations of $Q$ over an infinite sequence of meshes,
\begin{align}
\begin{split}
\mathbb{E}[\bar{Q}(\mathbf{y}; \tol)] &= \int_{\mathbb{R}^s} \bar{Q}(\mathbf{y}; \tol) \varphi(\mathbf{y}) d\mathbf{y}\\
&= \sum_{k = 0}^{\infty} \int_{Y_k, \tol} \bar{Q}_k(\mathbf{y}) \varphi(\mathbf{y}) d\mathbf{y},
\end{split}
\end{align}
which is a discrete approximation of the true integral
\begin{align}
\begin{split}
\mathbb{E}[Q(\mathbf{y})] &= \int_{\mathbb{R}^s} Q(\mathbf{y}) \varphi(\mathbf{y}) d\mathbf{y}.
\end{split}
\end{align}
The resulting bias satisfies:
\begin{equation}
	\abs{\mathbb{E}[\bar{Q}(\mathbf{y}, \tol)] - \mathbb{E}[Q(\mathbf{y})]} \leq \tol_{\textrm{bias}}. 
\end{equation}

In the following, we first analyze the single-level case, i.e., we consider a fixed tolerance, $\tol$ and suppress the $\tol$ dependence in the notation. Taking into account the expression for $a$ in~\eqref{eq:a_def}, the integrals w.r.t ${y}_1$ between discontinuity points admits an exact Gaussian CDF $\Phi$ representation for each fixed value of $\mathbf{y}_{-{1}}$, which allows us to bypass the need for quadrature approximation. 

The positions of the discontinuity points in ${y}_1$, denoted as $s_k$ for $k \in \mathbb{N}_0$, satisfy
\begin{align}
e^{-s_k\frac{p}{p+d}} \sum_K \tilde{\rho}_{k, K} h_{k, K}^{p+d} = \mathrm{TOL}_{bias} \frac{\int_{\mathcal{D}} {\tilde{\rho} }^{\frac{d}{p+d}} dx}{\int_{\mathcal{D}} \mathbb{E} \left[ \tilde{\rho}^{ \frac{d}{p+d}} \right] }, \quad  {k\in \mathbb{N}_0},
\label{eq:jump_mesh2}
\end{align}
where for each fixed $\mathbf{y}_{-1}$ we denote $\tilde{\rho} = \rho(0;\mathbf{y}_{-{1}})$, $\tilde{\rho}_{k, K} = \rho_{k, K}(0;\mathbf{y_{-{1}}})$, the error estimate with ${y}_1 \equiv 0$ on mesh $k$. The conditional expectation can be written as,
\begin{align}
	\label{eq:summation_Q}
\begin{split}
\mathbb{E}[\bar{Q}(\mathbf{y}) \vert \mathbf{y}_{-1}] &= \sum_{k = 0}^{\infty} \int _{s_k}^{s_{k-1}} \bar{Q}_{k} ({y}_1; \mathbf{y}_{-1})  \varphi({y}_1)d {y}_1\\
&= \sum_{k = 0}^{\infty} \int _{s_k}^{s_{k-1}} e^{-{y}_{1}} \tilde{Q}_{k} ( \mathbf{y}_{-1}) \varphi({y}_1)d {y}_1\\
&=  \sum_{k = 0}^{\infty} \tilde{Q}_{k} (\mathbf{y}_{-1}) \cdot \left( \Phi(s_{k - 1} + 1) - \Phi(s_k + 1) \right) \cdot \exp\left( \frac{1}{2} \right),
\end{split}
\end{align}
where we define $s_{-1} \coloneqq +\infty$. For the simplicity of the notation we also define
\begin{equation}
	\mu_k := \Phi(s_{k - 1} + 1) - \Phi(s_k + 1).
	\label{eq:mu_k}
\end{equation}

\section{Summation by parts}
\label{sec:summation_by_parts}

Notice that in the formulation~\eqref{eq:summation_Q}, $\mu_k$ does not converge to 0 uniformly w.r.t. $\tol$. This motivates us to apply summation by parts, a technique also discussed and analyzed in~\cite{vihola2018unbiased,glynn1992asymptotic,Rhee2015} to re-express the $\E{\bar{Q} | \mathbf{y}_{-1}}$:
\begin{equation}
\label{eq:summation_by_parts}
\begin{split}
\E{{\bar{Q}} | \mathbf{y}_{-1}} = \sum_{k=0}^{\infty} \tilde{Q}_k(\mathbf{y}_{-1}) \mu_k(\mathbf{y}_{-1}) &= \sum_{k=0}^{\infty} (\tilde{Q}_k(\mathbf{y}_{-1}) - \tilde{Q}_{k-1}(\mathbf{y}_{-1})) \left(\sum_{j=k}^{\infty} \mu_j (\mathbf{y}_{-1}) \right)\\
&:= \sum_{k=0}^{\infty} \Delta \tilde{Q}_{k} (\mathbf{y}_{-1}) \check{\mu}_k (\mathbf{y}_{-1}),
\end{split}
\end{equation}
for a fixed $\tol$. In this new expression~\eqref{eq:summation_by_parts}, $\check{\mu}_k (\mathbf{y}_{-1})$ monotonically decreases to 0 as $k \to \infty$, for a fixed $\tol$, which differs from $\mu_k$ in~\eqref{eq:summation_Q}. Notice that for a fixed $k$, $\check{\mu}_k \to 1$ as $\tol \to 0$. 

We consider a special case of the random field:
\begin{align}
	\label{eq:two_term_coefficient}
	a(\mathbf{x};\mathbf{y}) = \exp( \mathbf{y}_1 + \cos(\pi  \mathbf{x}_1) \sin(\pi \mathbf{x}_2) \mathbf{y}_2 ),
\end{align}
where ${y}_1, y_2 \sim \mathcal{N}(0,1)$ are independent. This choice of the coefficient is used to demonstrate the piecewise integration and summation-by-parts approaches discussed above. \revAE{The diagnostic study in this appendix uses a deliberately slower tolerance sequence, $\tol_\ell = 2^{-\ell - 2}$, so that component-level means and variances can be resolved over a wider range of~$\tol$; this is an independent single-level experiment on the two-term coefficient~\eqref{eq:two_term_coefficient} and should not be conflated with the multilevel runs of Section~\ref{sec:numex} on the Mat\'ern-like coefficient, which use $\tol_\ell = 4^{-\ell - 2}$.}

\begin{figure}[ht]
	\centering
	\includegraphics[width = 0.495\textwidth]{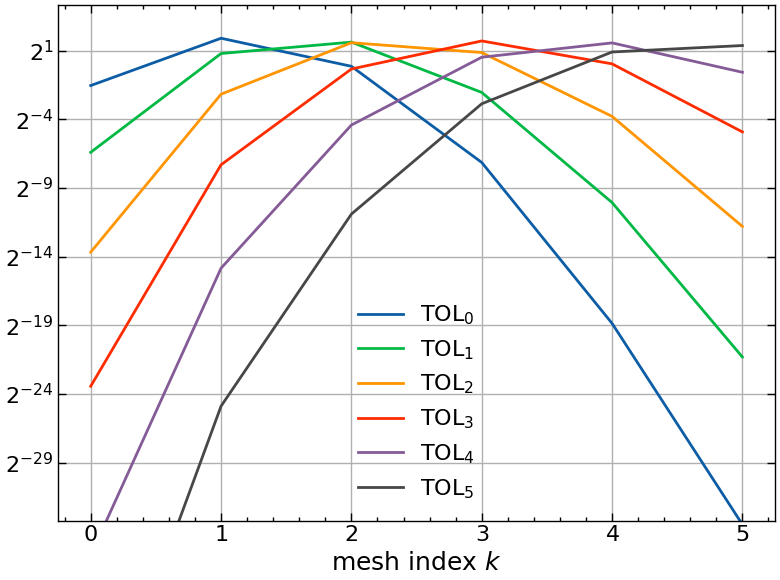}
	\includegraphics[width = 0.495\textwidth]{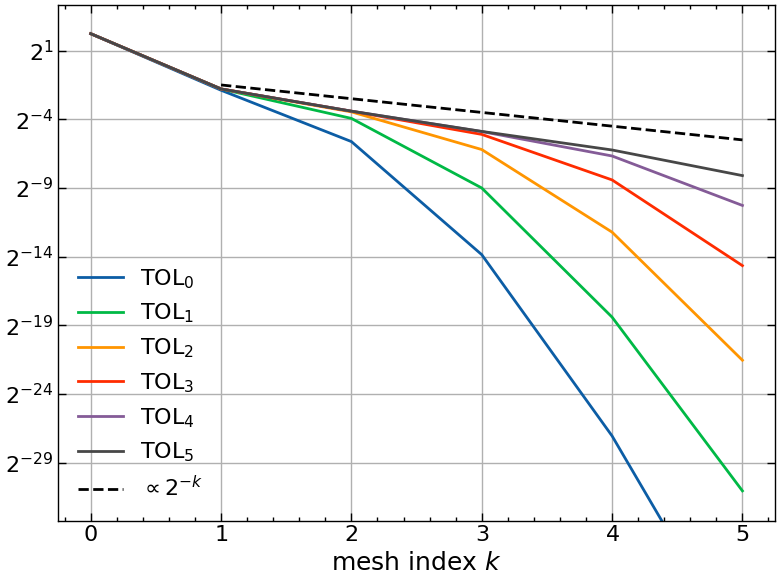}\\
	\caption{Example 1, Slit Domain: Expectation $\mathbb{E}\left[ \mu_k \tilde{Q}_k \right]$(left) and $\mathbb{E}\left[ \check{\mu}_k \Delta \tilde{Q}_k \right]$(right) for $k = 0, 1, 2, \dotsc, 5$ with $\tol_\ell = 2^{-\ell - 2}$, $\ell = 0, 1, 2, \dotsc, 5$ on Uniform Meshes. }
	\label{fig:single_level_mean_k}	
\end{figure}

Figure~\ref{fig:single_level_mean_k} displays $\mathbb{E}\left[ \mu_k \tilde{Q}_k \right]$ and $\mathbb{E}\left[ \check{\mu}_k \Delta \tilde{Q}_k \right]$. Notice that the dominating component index, $\arg \max_k \mathbb{E}\left[ \mu_k \tilde{Q}_k \right]$ increases as $\tol$ decreases, while $\arg \max_k \mathbb{E}\left[ \check{\mu}_k \Delta \tilde{Q}_k \right]$ remains 0, invariant of $\tol$. This behavior showcases the advantage of the summation-by-parts formulation, where $\check{\mu}_k$ converges uniformly to 0 as $\tol \to 0$. 

\begin{figure}[ht]
	\centering
	\includegraphics[width = 0.495\textwidth]{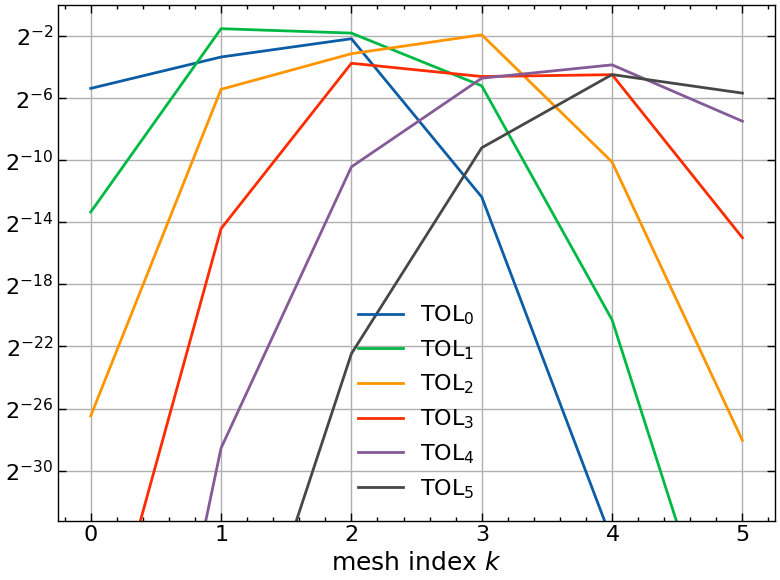}
	\includegraphics[width = 0.495\textwidth]{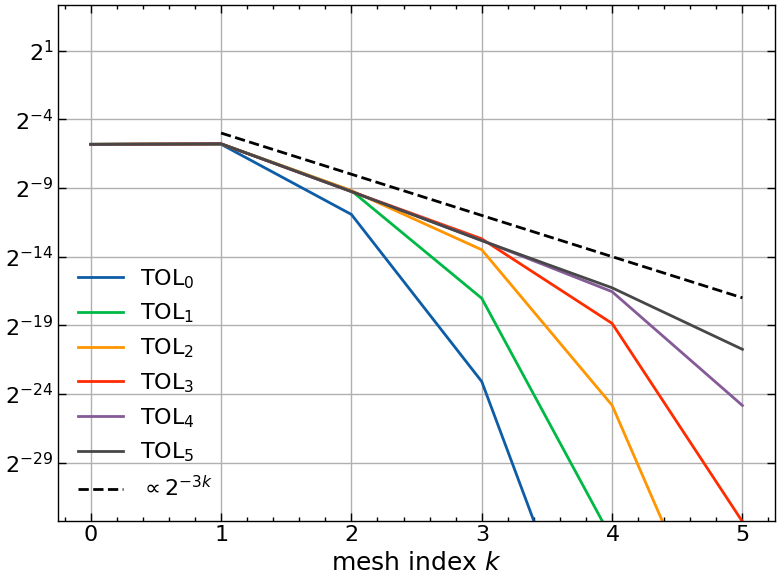}\\
	\caption{Example 1, Slit Domain: Variance $\var{\mu_k \tilde{Q}_k}$(left) and $\var{\check{\mu}_k \Delta \tilde{Q}_k }$(right) for $k = 0, 1, 2, \dotsc, 5$ with $\tol_\ell = 2^{-\ell - 2}$, $\ell = 0, 1, 2, \dotsc, 5$ on Uniform Meshes. }
	\label{fig:single_level_var_k_mc_qmc}
\end{figure}

\begin{figure}[ht]
	\centering
	\includegraphics[width = 0.495\textwidth]{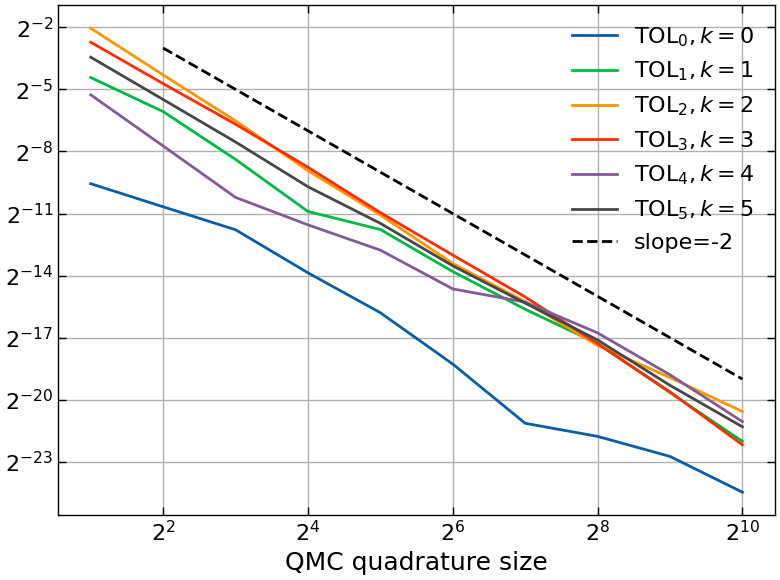}\\
	\includegraphics[width = 0.495\textwidth]{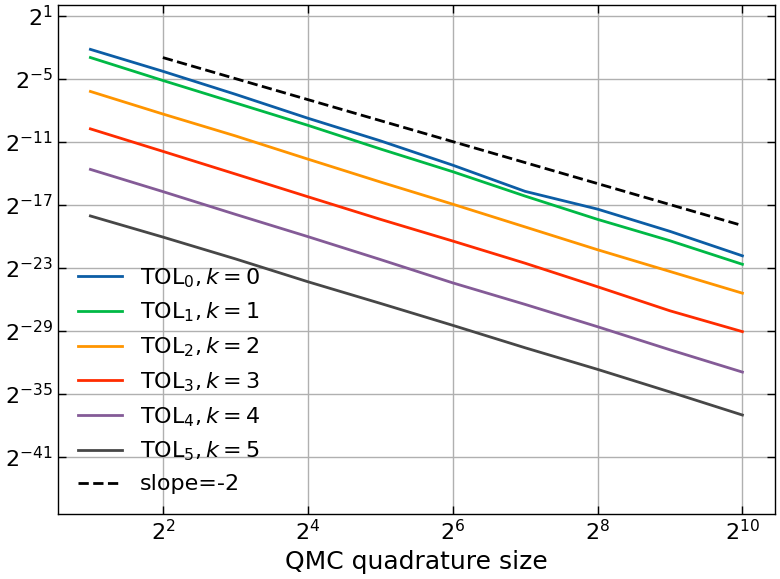}
	\includegraphics[width = 0.495\textwidth]{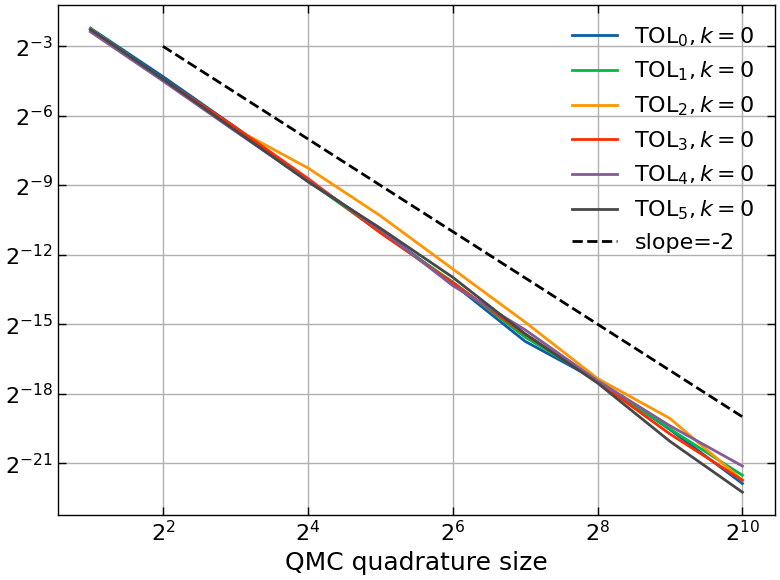}\\
	\caption{Example 1: Variance $\textrm{Var} \left[{I_{N_k} (\mu_k \tilde{Q}_{k}; \bm{\Delta}) }\right] $(top) and $\textrm{Var} \left[{I_{N_k} (\check{\mu}_k \Delta \tilde{Q}_{k}; \bm{\Delta}) }\right] $(bottom) for $k = 0, 1, 2, \dotsc, 5$ with $\tol_\ell = 2^{-\ell - 2}$, $\ell = 0, 1, 2, \dotsc, 5$, with 64 random shifts. }
	\label{fig:single_level_var_k_qmc}
\end{figure}

Figure~\ref{fig:single_level_var_k_mc_qmc} plots $\var{\mu_k \tilde{Q}_k}$ and $\var{\check{\mu}_k \Delta \tilde{Q}_k}$. The variance exhibits similar behaviors as those observed in the mean. For the RQMC method, Figure~\ref{fig:single_level_var_k_qmc} plots the decay of $\textrm{Var} \left[{I_{N_k} (\mu_k \tilde{Q}_{k}; \bm{\Delta}) }\right]$ and $\textrm{Var} \left[{I_{N_k} (\check{\mu}_k \Delta \tilde{Q}_{k}; \bm{\Delta}) }\right]$ with respect to some $k$. The slope of all curves is close to -2, consistent with the convergence rates derived in~\cite{liu2023nonasymptotic}. The above plot compares $\textrm{Var}  \left[{I_{N_k} (\mu_k \tilde{Q}_{k}; \bm{\Delta}) }\right]$ against $N_k$ with $k = \ell$ for each $\ell = 0,1,2, \dotsc, 5$. Although the rates of convergence are similar, for a given $N_k$, there is no clear pattern of $\textrm{Var}  \left[{I_{N_{k }} (\mu_{k } \tilde{Q}_{k }; \bm{\Delta}) }\right]$ when $\tol$ decreases. 

In the bottom plots of Figure~\ref{fig:single_level_var_k_qmc}, we show the decay of $\textrm{Var} \left[{I_{N_k} (\check{\mu}_k \Delta \tilde{Q}_{k}; \bm{\Delta}) }\right]$ against $N_k$ for with $k = \ell$ (left) and $k = 0$ (right). In this case, the variance for $k = \ell$, $\textrm{Var} \left[{I_{N_k} (\check{\mu}_k \Delta \tilde{Q}_{k}; \bm{\Delta}) }\right]$, decreases as $\tol$ decreases, for a same $N_k$, while the dominating component, $\textrm{Var} \left[{I_{N_0} (\check{\mu}_0 \Delta \tilde{Q}_{0}; \bm{\Delta}) }\right]$, remains nearly constant across different $\tol$.

\section{Importance Sampling}
\label{sec:mlqmc_is_appendix}
In this section we revisit the optimization problem to find out the optimal mesh size function $h^*$ within the framework of IS. 
{\paragraph{\revAE{Constrained optimization problem}}
The optimization problem adapted from the previous work~\cite{adMLMC_our} is stated as follows: Find the optimizer $h^{\ast}(\mathbf{y})$,
\begin{align}
\label{eq:optimization_problem_1}
\begin{split}
minimize \quad & \mathbb{E}_{\varphi_{\bm{\alpha}}} \left[ \int_{\mathcal{D}} h(\mathbf{y})^{-d} \right], \\
s.t. \quad &  \mathbb{E}_{\varphi_{\bm{\alpha}}} \left[ \int_{\mathcal{D}} \rho h^p(\mathbf{y}) \ell (\mathbf{y}) \right] \leq  \textrm{TOL}, \\
\end{split}
\end{align}
where $\ell({\mathbf{y}}) = \frac{\varphi(\mathbf{y})}{\varphi_{\bm{\alpha}(\mathbf{y})}}$. Notice that, compared to the settings in~\cite{adMLMC_our}, the integration measure changes from $\varphi$ to $\varphi_{\bm{\alpha}}$, the objective function remains the same as in the work~\cite{adMLMC_our}. However, the constraint now accounts for the IS, since IS modifies the integrand, thereby affecting the error estimate. The optimizer $h^{\ast}$ w.r.t $\mathbf{y}$ in~\eqref{eq:optimization_problem_1} is given by,
\begin{align}
h^{\ast}(\mathbf{x};\mathbf{y}) =  \frac{ \textrm{TOL}^{1/p}}{\left(\int_{\mathcal{D}}  \mathbb{E} \left[  (\ell \rho)^{\frac{d}{p+d}} \right] \right)^{1/p}} \ell(\mathbf{y})^{-\frac{1}{p+d}} \rho(\mathbf{x};\mathbf{y})^{-\frac{1}{p+d}}.
\label{eq:optimize_h_1}
\end{align}
With the optimal mesh function~\eqref{eq:optimize_h_1}, the error estimate satisfies the following equation:
\begin{align}
\label{eq:error_estimate_is_criteria}
\int_{\mathcal{D}} \rho(\mathbf{x};\mathbf{y}) {h^{\ast}}^p(\mathbf{x};\mathbf{y}) d\mathbf{x} =    \frac{ \textrm{TOL} }{\left(\int_{\mathcal{D}}  \mathbb{E} \left[ (\ell(\mathbf{y}) \rho(\mathbf{x};\mathbf{y}))^{\frac{d}{p+d}}  \right] d\mathbf{x} \right) } \int_{\mathcal{D}} \rho(\mathbf{x};\mathbf{y})^{\frac{d}{p+d}}  \ell(\mathbf{y})^{-\frac{p}{p+d}} d\mathbf{x}.
\end{align}
Notice that the error estimate~\eqref{eq:error_estimate_is_criteria} satisfies a different criterion than that from~\cite{adMLMC_our}. 
}



\bibliographystyle{siam}
\bibliography{bibliography_general, bibliography_MLQMC, bibliography_QMC_theory, bibliography_QMC_finance}

\end{document}